\documentclass{article}

\usepackage{algorithmic,amsfonts,amsmath,amsthm,tikz,float,url}

\usepackage[colorlinks=true,citecolor=blue]{hyperref}
\usepackage[capitalize]{cleveref}
\usepackage{caption} 

\usepackage{siunitx} 
\usepackage{multirow}
\usepackage{booktabs} 

\newtheorem{algorithm}{Algorithm}
\newtheorem{theorem}{Theorem}[section]
\newtheorem{lemma}[theorem]{Lemma}
\newtheorem{corollary}[theorem]{Corollary}

\newtheorem{remark}[theorem]{Remark}

\newtheorem{proposition}[theorem]{Proposition}
\newtheorem{conjecture}[theorem]{Conjecture}
\newtheorem{example}[theorem]{Example}
\newcommand{\RN}{\textnormal{RN}}

\newcommand{\ulp}{\textnormal{ulp}}
\newcommand{\emin}{{e_\mathrm{min}}}
\newcommand{\emax}{{e_\mathrm{max}}}

\newcommand{\umax}{{u_\mathrm{max}}}

\usepackage{fancyvrb}
\fvset{fontshape=normal}
\fvset{listparameters={\setlength{\topsep}{0pt}\setlength{\partopsep}{0pt}}} 

\floatstyle{ruled}
\newfloat{algorithm}{loa}{List of Algorithms}
\floatname{algorithm}{Algorithm}

\allowdisplaybreaks 

\title{Effective Quadratic Error Bounds for
Floating-Point Algorithms Computing the Hypotenuse Function}
\author{Jean-Michel Muller \and Bruno Salvy}

\begin{document}
\maketitle

\begin{abstract}
We provide tools to help automate the error analysis of algorithms 
that evaluate simple
functions over the floating-point numbers. The aim is to obtain
tight relative
error bounds for these algorithms,
expressed as a function of the unit round-off. Due to the discrete nature
of the set of floating-point numbers, the largest errors are often intrinsically ``arithmetic'' in the sense that their appearance may depend on specific bit patterns in the binary representations of intermediate variables, which may be present only for some precisions. 
We focus on \emph{generic} (i.e., parameterized by the precision) and \emph{analytic} over-estimations that still capture
the correlations between the errors made at each step of the
algorithms. Using methods from computer algebra, which we adapt to  the particular structure of the polynomial systems that encode the
errors, we obtain bounds with a linear term in the unit round-off that
is sharp in many
cases. An explicit quadratic bound is given, rather than the $O
()$-estimate that is more common in this area. This is particularly important when using
 low precision formats, which are increasingly common in modern
processors. Using this approach, we compare five algorithms for computing the hypotenuse function, ranging from elementary to
quite challenging.
\end{abstract}
\section{Introduction}

\subsection{Motivation}
Floating-Point (FP) arithmetic is ubiquitous in numerical computing: almost all High-Performance Computing applications rely on it. However, FP arithmetic is inherently inexact. The impact of  individual rounding errors on the final result of a computation is very often small but can sometimes be catastrophic (see e.g.,~\cite[Chapter~1]{AltmannGillMcDonald2004}, \cite{GAO1992}). It is therefore important to be able to obtain bounds on the error that can occur when running a given numerical program.  

Assuming all variables remain in the so-called ``normal domain'', the
relative error of each individual arithmetic operation can be bounded
by a constant $u$, that depends only on the floating-point format
being used (more detail is given in Section~\ref{subsec:FParith}). 
Consequently, the natural and most commonly adopted way to obtain
error bounds is to ``propagate'' these individual bounds iteratively.
For example, the following very simple program computes the product of three FP numbers $a$, $b$, and $c$ by two successive multiplications:
\begin{verbatim}
d = a * b;
f = d * c
\end{verbatim}
The computed value of \texttt{f} satisfies
\[
abc(1-u)^2 \leq f \leq abc(1+u)^2.
\]
This simple approach has proved fruitful since the pioneering work of
Wilkinson~\cite{Wilkinson1960}. It has made it possible to obtain
numerous results~\cite{Hig02}. Unfortunately, it suffers from two  drawbacks:
\begin{enumerate}
    \item Very often, it leads to a large overestimation of the real
    maximum error. First, because some operations are errorless:
    examples are the subtraction of FP numbers that are close to each
    other (Lemma~\ref{lemma-sterbenz}), or multiplications by powers
    of $2$. Second, and probably more important, the bound $u$ on
    the relative error of an operation is sharp only if its result is
    very close to (and above) a power of $2$ (see \cref{subsec:FParith}). Even for
    fairly simple programs, this cannot happen for all intermediate variables because they cannot realistically be viewed as ``independent'' (just consider computing $3xy$ by first computing $xy$ and then multiplying the resulting value by $3$: you cannot have $xy$ and $3xy$ very close to a power of $2$ at the same time);
     \item When the size of the analyzed program becomes large (more than a few operations), the expressions obtained by propagating the individual error bounds become too large and complex to be easily manipulated by paper-and-pencil calculations, with several consequences. The first one is that  the proofs of the theorems that give the bounds are long, tedious, and thus error-prone, with the
    unpleasant consequence that one is never quite sure of their
    correctness because few people read them in detail. The second
    consequence is that to avoid this complexity, it is tempting to ``simplify'' the intermediate bounds at each step, so that they become looser.
\end{enumerate}
Our goal is to alleviate these drawbacks by (partly) automating the computation of error bounds, using modern computer algebra tools. We believe that this approach would be beneficial in a number of ways:

\emph{Reusability:} when designing and analyzing numerical algorithms,
one often wants to compare slightly different variants. ``Re-playing'' the automatic computation of the error bound with a new variant is easily done.

\emph{Tighter bounds:} modern computer algebra tools can readily manipulate complex expressions, so that intermediate simplifications are much less often needed, with the result that the final bounds are often tighter than those obtained from paper-and-pencil calculations.

\emph{Trust:} automatically generating
    bounds, together with their proofs, which are correct by
    construction would greatly reduce the likelihood of an error  remaining unnoticed. Faced with a
    similar problem, Muller and Rideau
    use formal proofs~\cite{MullerRideau2022}.
    This is also the case for Gappa~\cite{DaumasMelquiond2010} and \texttt{Real2Float}~%
    \cite{MagronConstantinidesDonaldson2017} mentioned below. These approaches complement our own, and a natural next step
    will be to have the computer algebra tool generate a certificate
    to be validated by the formal proof system. At this stage however,
    we are not sufficiently confident in our implementation to trust
    its results blindly. For this reason, we give both the results
    given by our implementation and a human-readable proof of them, or
    in some cases  weaker bounds that those found automatically.

\subsection{Basics of Floating-Point Arithmetic}
\label{subsec:FParith}
We give the definitions that are relevant for this study. We refer to
surveys and books for more
information on floating-point arithmetic~\cite{BoldoEtAl2023,Gol91,MullerEtAl2018,Overton2001}.
Throughout this article, we assume a radix-$2$  FP system 
parameterized by its \emph{precision} $p$  and its  \emph{extremal exponents} $\emin$ and $\emax$. That is, an FP number is of the form
\begin{equation}
\label{def-fp-1}
x = M \times 2^{e-p+1},
\end{equation}
where $M$ and $e$ are integers satisfying 
\begin{equation}
\label{def-fp-2-3}
|M| \leq 2^p-1 \qquad\text{and}\qquad
\emin \leq e \leq \emax.
\end{equation}
The FP number $x$ is said \emph{normal} if $|x| \geq 2^\emin$, and 
\emph{subnormal} otherwise. 
The largest finite FP number is
$
\Omega = (2^p-1) \cdot 2^{\emax-p+1},
$
and the \emph{normal domain} is the range
$2^\emin\le|t|\le\Omega.$

As the exact sum, product and quotient of two FP numbers (or the square root of a FP number) are not, in general, FP numbers, they must be \emph{rounded}. In the following, we assume that the rounding function is \emph{round-to-nearest}, noted $\RN$, which is the default\footnote{More precisely, the  default in the IEEE-754 Standard is
round-to-nearest \emph{ties-to-even}, see for instance~\cite{MullerEtAl2018} for more explanation.} in the IEEE 754 Standard on Floating-Point Arithmetic~\cite{IEEE754-2019}. That is, each time $a \star b$ (resp. $\sqrt{a}$) is computed, where $a$ and $b$ are FP numbers and $\star \in \{+,-,\times,\div\}$, what is returned is $\RN\left(a \star b\right)$ (resp. $\RN\left(\sqrt{a}\right)$). The absolute error due to rounding to nearest a real number $t$, $|t| \leq \Omega$ is bounded by half the distance between two consecutive FP numbers in  the neighborhood of $t$. That distance, called \emph{unit in the last place} (ulp) of $t$ is
\begin{equation}
\label{def-ulp}
\ulp(t) :=
\left\{
\begin{array}{lll}
2^{\emin-p+1} & \mathrm{if~} |t| < 2^{\emin+1} \\
2^{\left\lfloor \log_2 |t| \right\rfloor -p+1} & \mathrm{otherwise}.
\end{array}
\right.
\end{equation}
Assume that $t$ belongs to the normal domain. There exists $k \in \mathbb{Z}$, $k \geq \emin$ such that $t \in [2^k,2^{k+1})$. The \emph{absolute error} due to
    rounding $t$ is bounded by
    \begin{equation}
    \label{err-abs}
    \left|t - \RN(t)\right| \leq \frac{1}{2} \ulp(t) = 2^{k-p},
    \end{equation}
   and therefore, the
    \emph{relative error} due to rounding $t$ is bounded by
    \begin{equation}
    \label{rel-err-u}
   \left|\frac{t - \RN(t)}{t}\right| \leq  u,
   \end{equation}
   where $u := 2^{-p}$ is the \emph{unit round-off}.
   
\begin{figure}
\centerline{\includegraphics[width=0.4\textwidth]{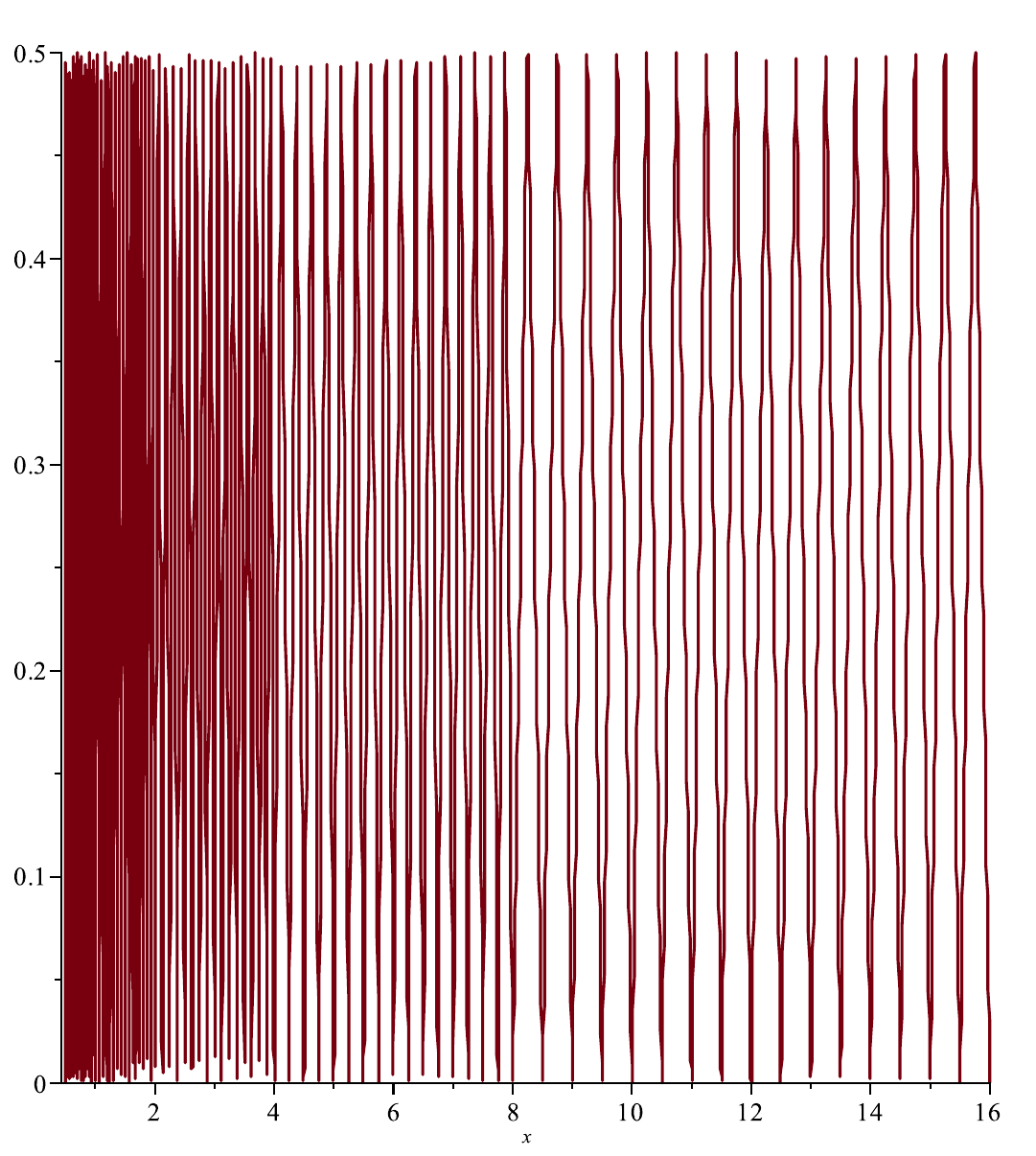}
\hfil\includegraphics[width=0.45\textwidth]{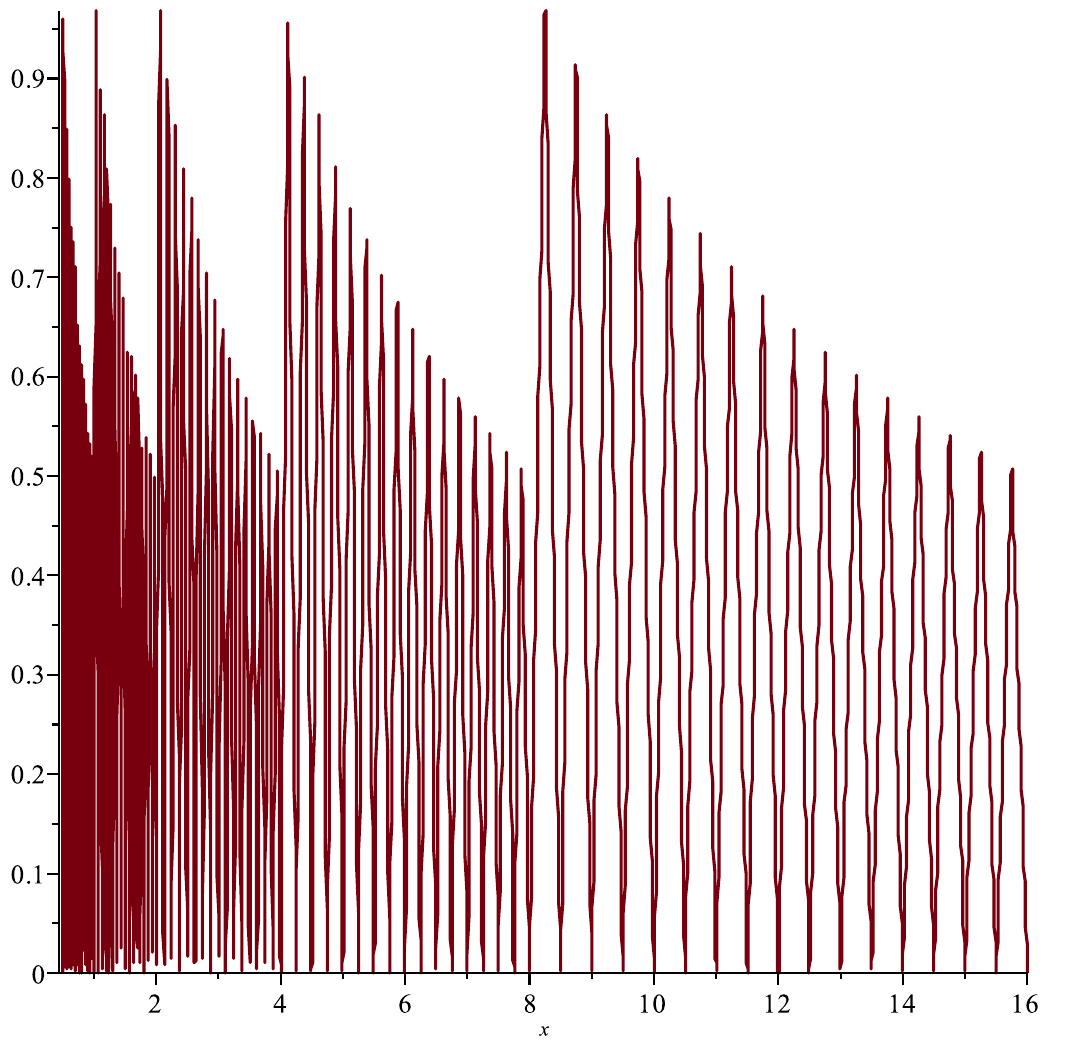}}
\caption{Left: absolute error (in ulps) of rounding to nearest $x \in [\frac{1}{2},16]$.
Right: relative error (in multiples of $u=2^{-p}$) of rounding to
nearest $x \in [\frac{1}{2},16]$. Both pictures assume a
binary floating-point system with $p=5$.}
\label{fig-err}
\end{figure}

It is obvious that  \cref{err-abs} $\Rightarrow$ \cref{rel-err-u}, but
the reverse is not true. Still, although it conveys less
information than \cref{err-abs},  \cref{rel-err-u} is used more often than \cref{err-abs} in both
paper-and-pencil and automated error analysis, because it is easier to
manipulate relative errors when analyzing long sequences of
operations.  However,  it is unlikely that near-optimal error bounds can be obtained using only \cref{rel-err-u}, except for very short and simple algorithms. Figure~\ref{fig-err} shows the absolute and relative errors resulting from rounding a real number $t \in [1/2,16]$. The absolute error bound $0.5 \ulp$ is approached in the immediate vicinity of all real numbers, while the relative error  bound $u$ is approached only slightly above powers of~$2$. As a consequence, obtaining a near-optimal bound using \cref{rel-err-u} solely requires the existence of input values for which almost all intermediate variables in the computation under consideration are slightly above a power of~$2$. While this may happen for very simple algorithms such as Algorithm~\ref{NaiveHypot} below, it is unlikely to happen for more complex algorithms.

\subsection{Which kind of error bounds?}
\label{sec:kindofbounds}
\paragraph{Generic and analytic bounds}
Our aim is to determine relative error bounds that we term as \emph{generic} and \emph{analytic} for algorithms that
typically contain at most a dozen or so arithmetic operations. We use these terms to denote that the bounds are expressed as a continuous function of  $u$ that bounds the error for all $u$ that are less than some specified $\umax > 0$. 
To be more specific, we are interested in obtaining \emph{generic
quadratic bounds}, i.e., generic analytic bounds of the form
\begin{equation}
   \label{desired-form}
   \alpha u+\beta u^2,\qquad u\le \umax.
\end{equation}
We look for ``best possible'' generic quadratic bounds: for $u \leq
\umax$ we want to minimize $\alpha$, and whenever possible, for that best $\alpha$ we want to minimize $\beta$. {In particular, we do not want to return \emph{approximate} bounds, such as those of the form ``$\alpha u + \mathcal{O}(u^2)$'' that are frequent in numerical analysis,\footnote{Such approximate bounds are of course very useful in many applications, but we mainly target the small algorithms that implement the ``basic building blocks'' of computing, for which guaranteeing upper bounds on the error is important.} because they do not
allow to \emph{guarantee} that the error will be less than some well-defined value.} To accomplish this, as said above,
utilizing solely the relative error model  (\ref{rel-err-u}) will
not always be adequate. Instead,  we often need to use the absolute-error model  (\ref{err-abs}), in conjunction with particular properties of the set of FP numbers such as Lemma~\ref{lemma-sterbenz} or \cref{lemma:jeannerod-rump} below.

Consider, for a given algorithm an (absolute or relative) error bound
$\mathcal{B}(u)$ and the (most likely unknown) worst case error $
\mathcal{W}(u)$. We call the bound $\mathcal{B}$
\begin{itemize}
    \item \emph{asymptotically optimal} if  $\mathcal{W}(u)/\mathcal{B}(u) \to 1$  as $u \to 0$;
    \item \emph{sharp} (for $u \leq \umax$) if $\mathcal W(u)\ge 0.999 \mathcal B(u)$ for some values
    of~$p$\footnote{The constant $0.999$ is of course arbitrary.}.
\end{itemize}

The exact motivation for a careful error analysis and the underlying hypotheses it uses may vary greatly from author to author, and this can influence the type of bound one looks for. As regards motivation for computing numerical error bounds, one can cite:
\begin{itemize}
    \item \emph{Choosing between different algorithms:} if two different
    algorithms are available to solve the same problem, one may wish to make an
    informed choice of the algorithm that has the best balance
    performance/accuracy. This requires  that the error bounds be sharp (but approximate bounds could do).
    \item \emph{Careful implementation optimization:}  There is a comprehensive range of FP formats, ranging from $16$-bit ``half-precision'' formats to the $128$-bit ``quad-precision'' binary128 format. Since numbers represented in the small formats are processed faster, there is a temptation to use the small formats whenever possible~\cite{abdelfattah2020survey}. As the small formats have significantly larger individual rounding errors, this requires sharp error bounds.
    
    \item \emph{Certainty:} floating-point arithmetic is also used in critical numerical software (e.g., software embedded in transportation systems). For critical applications, it is sometimes essential to be sure that the numerical error does not exceed a certain limit. In such cases, the sharpness of the error bounds is not always needed, but certainty is paramount: approximate bounds are to be avoided. 
\end{itemize}

\paragraph{Genericity versus optimality}

The \emph{genericity} of the error
bounds avoids having to repeat the analysis for all possible FP formats. However, looking for generic bounds may mean giving up on optimality. The ``best possible'' generic analytic bound will sometimes not be sharp, because rounding errors  are inherently ``arithmetic'': they may depend on specific bit patterns in the binary representations of intermediate variables, which may be present only for some  values of $p$. Consider for example the computation of $x^2-2$ using a  multiplication followed by a subtraction, i.e., what is actually computed is
$
\RN(\RN(x^2)-2)).
$
Under the (generally accepted) conjecture~\cite{BaileyCrandall2001}
that $\sqrt{2}$ is a $2$-normal number\footnote{There is an
unfortunate conflict of scientific terminology between the fields of
computer arithmetic and number theory here: the word ``normal'' used
here has no correlation with its usage in other parts of the article.}
(which implies than any given bit string appears infinitely many times
in its binary representation), the behaviour of the algorithm is
very different for~$x=\RN(\sqrt{2})$ and the other values. The
following properties hold:
\begin{itemize}
    \item[(P1)] in precision-$p$ arithmetic, if $x \neq \RN(\sqrt{2})$ then the relative error of the computation is less than or equal to $\sqrt{2}/2 \approx 0.707$;
    \item[(P2)] for infinitely many values of $p$, the largest relative error of the computation (hence attained for $x = \RN(\sqrt{2})$) is 1. This is found by choosing $p$ such that just after the first $p$ bits of the binary expansion of $\sqrt{2}$ there is either the bit-string $000$ or the bit string $111$ (the normality conjecture implies that there are infinitely many such $p$). In such cases, for $x = \RN(\sqrt{2})$, $\RN(x^2) = 2$ so that the computed result is zero;
    \item[(P3)] for infinitely many values of $p$, the relative error of the computation for $x = \RN(\sqrt{2})$ is less than $128/(112\sqrt{2}-49u)$, which is less than $0.876$ as soon as $u \leq 1/4$ (i.e., $p \geq 2$). This is found by choosing $p$ such that just after the first $p$ bits of the binary expansion of $\sqrt{2}$ there is either the bit-string $10000$ or the bit string $01111$.
\end{itemize}

Let $\mathcal{B}(u)$ be a generic analytic bound. (P2) implies that there are values of $u$ arbitrarily close to $0$ for which the largest error is $1$. For these values we must have $\mathcal{B}(u) \geq 1$ and since  $\mathcal{B}$ is continuous by hypothesis,  $\mathcal{B}(0) \geq 1$ so that for any $\epsilon > 0$ there exists $u_0 > 0$ such that for $u < u_0$, we have $\mathcal{B}(u) \geq 1-\epsilon$. But (P1) and (P3) imply that there are infinitely many values of $p$ for which the largest relative error is less than $0.876$, implying that $\mathcal{B}(u)$ is not asymptotically optimal (and not even sharp).

\subsection{Recent results}
\label{sec:recentresults}
\paragraph{Recent results in computer arithmetic}
Rump revisited the classical error bounds for recursive summation and
dot product~\cite{Rump:2012cu},  showing that the usual ``$\mathcal{O}(u^2)$'' terms which appear in the literature are not necessary. This prowess led to similar improvements for other problems such as summation in arbitrary order, polynomial evaluation, powers and iterated products, LU factorization, etc~\cite{Rump2019}.
In this vein, Jeannerod and Rump proved optimal bounds for square-root
and division~\cite{JeannerodRump2018}, recalled in \cref{lemma:jeannerod-rump}
below.
These bounds offer precise control over the accuracy of
fundamental operations. They are obtained through a delicate analysis of
the discrete structure of the set of the FP numbers, which we do
not attempt to automate here.

\paragraph{Recent work on automatic analysis}
Various tools have been proposed for computing error bounds on the
result of numerical programs. Noteworthy tools include \texttt{Fluctuat}~%
\cite{Goubault2007} (based on abstract interpretation), \texttt{FPTaylor} and
\texttt{SATIRE}~\cite{Solovyev2015,DasEtAl2020} (based on Taylor forms), \texttt{Gappa}~%
\cite{DaumasMelquiond2010} (based on interval arithmetic and forward
error analysis), \texttt{VCFloat2}~\cite{AppelKellison2024} 
(based on interval arithmetic and a Coq tactic that recursively decomposes expressions), and \texttt{Real2Float}~%
\cite{MagronConstantinidesDonaldson2017} (based on semidefinite
programming and the use of  (\ref{rel-err-u})). They are very
efficient in
their
respective
domains
(for instance Fluctuat and SATIRE can address rather large programs,
Gappa and Real2Float can generate formal proofs or formally provable
certificates, thence giving very good confidence
on the obtained results). In general, they are somehow limited to a
given precision (e.g., double-precision) and cannot return a bound of
the form~(\ref{desired-form}) valid with any $u \leq \umax$. Moreover,
the ability of handling large programs comes at the expense of
bounds that may significantly exceed the optimal ones.

\subsection{Contribution}

Our approach is complementary to the works mentioned above.
We
provide generic analytic bounds
parameterized by the unit round-off~$u$. They are computed in such a way that some of
the correlations between the errors in the intermediate steps are
taken into account. This often yields tighter bounds than those
provided by previous tools, 
at the expense of heavy computations. This limits our approach to the small programs that implement functions considered
as ``basic building blocks'' of numerical computation, such as the hypotenuse considered in this
article, for which it makes sense to spend much computing time in order to derive a sharp bound that will be reused many times, in
many applications.

\paragraph{Running example: the hypotenuse function}
We illustrate our approach with a gallery of algorithms for the computation of  $\sqrt{x^2+y^2}$, presented in Section~\ref{sec:algos-hypot}. A basic example is provided by the following naive program:
\begin{verbatim}
sx = x*x;
sy = y*y;
sigma = sx + sy;
rho1 = sqrt(sigma);
\end{verbatim}

Assuming that the compiler does not change the sequence of
instructions of the program, and that the arithmetic operations and the square root are rounded to nearest, what is actually computed is expressed by
Algorithm~\ref{NaiveHypot} below. 

\begin{algorithm}[H]
 \caption{The naive algorithm for the \texttt{Hypot} function. It takes $4$ FP operations, and approximates $\sqrt{x^2+y^2}$ with relative error better than $2u$.}
\label{NaiveHypot}
  \begin{algorithmic}[1] 
    \STATE{$s_x \gets \RN(x^2)$}
    \STATE{$s_y \gets \RN(y^2)$}
    \STATE{$\sigma \gets \RN(s_x + s_y)$}
    \STATE{$\rho \gets \RN(\sqrt{\sigma})$}
  \end{algorithmic}
\end{algorithm}

Assuming no underflow or overflow occurs, and just using 
\cref{rel-err-u} on each statement yields the relative error bound
$2u+u^2$.
This is the classical $2u + \mathcal{O}(u^2)$ error bound of the
literature~\cite{HFT94}.  Jeannerod and Rump improved that
bound by showing that the relative error is actually less than 
$2u$~\cite{JeannerodRump2018}. Both bounds are {asymptotically optimal}~\cite{JeannerodMullerPlet2017}. 
As a first example of our approach, we show in \cref{sec:warm-up} that
the error is bounded by 
\begin{equation}\label{eq:bound-algo1}
2u-\frac{8}5(9-4\sqrt6)u^2<2u-\frac54u^2,\qquad u\le 1/4.
\end{equation}

\paragraph{Main tool: Polynomial Optimization}
By the very nature of the FP numbers, the error is a
discontinuous function of $u$ (and $u = 2^{-p}$ itself only takes discrete values). However, the sharp error bounds that are known
for the basic operations are all continuous functions  of $u$ and, even better from the
computer algebra point of view, algebraic functions  (being obtained by
combinations
of $+,\times,\div,\sqrt{}$). Thus an algorithm can be seen as the
description of a
semi-algebraic set: a subset of~$\mathbb R^n$ defined by polynomial
equalities and inequalities. In the specific example of
Algorithm~\ref{NaiveHypot}, the equations
are
\begin{equation}\label{eq-algo1}
\begin{split}
s_x&=x^2(1+u\epsilon_{s_x}),\quad s_y=y^2(1+u\epsilon_{s_y}),\\
\sigma&=(s_x+s_y)(1+u\epsilon_\sigma),\quad\alpha^2=\sigma,\quad
\rho=\alpha(1+u\epsilon_\rho),
\end{split}
\end{equation}
with each $\epsilon_i$ bounded in absolute value by~$1/(1+u)$. In
general, tight bounds on these quantities~$\epsilon_i$ are obtained
by a step-by-step analysis of the program: in some cases they 
follow from a refinement of the relative error bound \cref{rel-err-u};
in other cases
knowledge of an interval containing the value allows the use of the
tighter absolute error bound \cref{err-abs}.
In our example, this system of equalities and inequalities
defines a small tube-like region of dimension~7 in~$\mathbb R^{12}$:
here, the
dimension~12 comes from one
variable
for the unit round-off~$u$, one
for each input~$x,y$, one for each of the variables of the algorithm 
($s_x,s_y,\sigma,\rho$), one for each $\epsilon_i$ and one for the
square-root~$\alpha$.
A computation reveals that the projection of this region on the $
(u,x,y,\rho)$-coordinate space is
simply the
region between the surfaces $\rho=\sqrt{x^2+y^2}/(1+u)^2$ and $\rho=
\sqrt{x^2+y^2}(1+2u)^2/(1+u)^2$. From there, it is easy to deduce the bound in
\cref{eq:bound-algo1}.

For less straightforward algorithms,
obtaining a bound becomes an optimization problem: find the
maximum and the minimum of $\rho/\sqrt{x^2+y^2}-1$ in the domain
defined by
the equations and inequalities obtained from analyzing the program. 
This is an instance of the \emph{polynomial optimization problem}, 
which consists in
computing the maximum or the minimum of a variable that satisfies a
given system of polynomial equations and inequalities. 
At this level of generality, this problem is
well understood in terms of algorithms and complexity, see e.g., 
\cite{NieRanestad2009,JeronimoPerrucciTsigaridas2013,BankGiustiHeintzEl-Din2014}, but also very
expensive. For large algorithms, one has no choice but to compute
approximations to the optimal values; this is the approach taken by
\texttt{Real2Float}~\cite{MagronConstantinidesDonaldson2017}. We focus here
on testing the limits of what can be computed exactly. 
Our approach exploits the specific
structure of the problem as much as we can and is described in \cref{app-ca}.

\paragraph{Prototype Implementation}
We have written a simple Maple implementation of the approach described in this
article\footnote{A Maple session  available on arXiv with this article
includes all its examples.}.
On the example above, its
behaviour is as follows:
\begin{Verbatim}
> Algo1:=[Input(x=0..2^16,y=0..2^16,_u=0..1/4), 
> s[x]=RN(x^2),s[y]=RN(y^2),sigma=RN(s[x]+s[y]),rho=RN(sqrt(sigma))]:
> BoundRoundingError(Algo1);
\end{Verbatim}
\[{\color{blue}2\textit{\_u}+\left(\frac{72}5 - \frac{32\sqrt{6}}5\right)\textit{\_u}^2}\]
The algorithm is first stated with ranges for its inputs (currently the program
only handles finite ranges, which is not much of a problem for the hypotenuse function, but might make the analysis somehow difficult for other functions). The name \verb|_u| is used for the unit round-off.
The procedure first makes several decisions
concerning intermediate rounding errors; next, it
computes a bound on the linear part of the error bound.
Finally, the procedure 
computes an upper
bound on the quadratic part once this linear part is fixed. 
This is the bound from \cref{eq:bound-algo1}. Thus in a way the input of our
code is an algorithm and its output is a theorem giving an upper bound on its
relative error. However, as mentioned above, since this is only a
prototype, in this article we do not
trust the output blindly. Instead, we use it as a statement of a theorem and
give a paper proof, obtained mostly from following the intermediate
steps of the code, so that the proof can be checked by a human reader.

\subsection{Structure of the article}

In \Cref{sec:intro-bounds}, we recall special cases where the
bounds \cref{err-abs} and \cref{rel-err-u} on the errors of individual
operations can be improved due to some structural properties of the
set of the FP numbers. Then, \Cref{app-ca} describes the computer
algebra algorithms we use for analyzing numerical programs. \Cref
{sec:algos-hypot} presents several algorithms that have been proposed
in the literature for evaluating the hypotenuse function. These
algorithms are analyzed in Sections~\ref{sec:warm-up} 
to~\ref{sec:analysis-Kahan}. We discuss these results in 
\cref{sec:discussion}. The most technical parts of the analyses are
deferred to Appendices \ref{app-proof-beebe} to \ref{sec:proof-kahan}.

\paragraph{Acknowledgements}  This work is partly supported by the
ANR-NuSCAP 20-CE-48-0014 project of the French \emph{Agence nationale
de la recherche} (ANR). We are grateful to Mohab Safey El Din for
giving us access to recent versions of his software, to Guillaume
Melquiond for his help with Gappa and to Ganesh Gopalakrishnan and
Tanmay Tirpankar for their help with Satire.

\section{Bounds for Floating-Point Operations}\label{sec:intro-bounds}

Some structural properties of the set of FP
numbers cannot be disregarded if tight
error bounds are desired. An illustration is given by Lemma~\ref{lemma-sterbenz}
below: if two FP numbers $a$ and $b$ are close enough, then computing $a-b$ results in no error. In such a case, utilizing  \cref{err-abs} or \cref{rel-err-u} to bound the error  would be a significant overestimation. A simpler example is multiplications by powers of two, which are exact operations as long as underflow and overflow are avoided.

\begin{lemma}[Sterbenz' Lemma~\cite{Ste74}]
\label{lemma-sterbenz}
If $a$ and $b$ are floating-point numbers satisfying $a/2 \leq b \leq 2a$ then $b-a$ is a floating-point number, which implies $\RN(b-a) = b-a$.
\end{lemma}

In a similar spirit, the following result of Boldo and 
Daumas~\cite[Thm.~5]{BolDau03a} is helpful to deal with the error in a square-root
computation.
 
 \begin{lemma}[Exact representation of the square root remainder \cite{BolDau03a}]
\label{lemma-boldodaumasSQRT}
In binary, precision-$p$, FP arithmetic with minimum exponent $\emin$, let $s = \RN(\sqrt{t})$, where $t$ is a FP number. The  term
$
t - s^2
$
is a FP number if and only if there exists a pair of integers $(m,e)$ (with $|m|
\leq 2^p-1$) such that $s = m \cdot 2^{e-p+1}$ and $2e \geq \emin+p-1$.
\end{lemma}

An almost immediate consequence of \cref{def-ulp} is the
following slight improvement of \cref{rel-err-u}.
\begin{lemma}[Dekker-Knuth's bound~\cite{Knu98}]
\label{lemma-knuth-dekker}
If $t\neq0$ is a real number in
the normal domain, then the relative error due to rounding satisfies
\begin{equation}
\label{rel-err-v}
\left|\frac{t-\RN(t)}{t}\right| \leq v\qquad\text{with}\quad v = 
\frac{u}{1+u}.
\end{equation}
\end{lemma}

This bound was obtained by Dekker in 1971~\cite{Dekker71} for the
error of some operations in binary FP arithmetic. It was given by
Knuth in its full generality~\cite{Knu98}. But it is seldom used:
most authors
use the very slightly looser (but simpler)  bound~(\ref{rel-err-u}).
In our context where the analyses are performed by a computer, we use
\cref{rel-err-v} when appropriate.

\paragraph{Relative error for specific operations}
Jeannerod and Rump~\cite{JeannerodRump2018} have showed that while the
bound (\ref{rel-err-v}) is optimal for addition, subtraction and
multiplication, it can be improved for division and square root.
For these two operations, they give the following optimal bounds, of
which we make
heavy use.
\begin{lemma}[Jeannerod-Rump bounds~\cite{JeannerodRump2018}]
\label{lemma:jeannerod-rump}
When the precision $p$ is at least~2, the relative
error of a (correctly rounded) square root is bounded by
\begin{equation}
\label{err-sqrt-jr}
1 - \frac{1}{\sqrt{1+2u}};
\end{equation}
the relative error of a division in binary FP arithmetic is bounded
by
\begin{equation}
\label{err-div-jr}
u-2u^2.
\end{equation}
\end{lemma}
The first bound is general; the second one holds only in base~2, which
is our setting in this article.

\section{Computer algebra algorithms for tight analysis of numerical
programs}\label{app-ca}
We focus here on the analysis of \emph{straight-line programs}
operating over floating-point numbers. We assume that the basic operations in these programs are $+$, $-$, $\times$, $\div$, $\sqrt{.}$ and the FMA operation (which evaluates an expression $ab+c$ with one final rounding only and is available on all recent processors). Hence, such a program is a sequence of
instructions, the $i$th one of which has one of the forms
\begin{equation}\label{eq:instruction}
 v_i:=w_j\times w_k+w_\ell\quad\text{or}\quad v_i:=w_j\div
 w_k\quad\text{or}\quad
v_i:=\sqrt{w_j},
\end{equation}
where  $w_j,w_k,w_\ell$ are either variables~$v_m$ with
$m<i$, or input variables, or constant floating-point numbers.
For the first of the three forms (which corresponds to the FMA operation), taking some of the variables to be~0 or~1 recovers addition and
multiplication. 
By convention, the
last instruction defines the result of the algorithm, whose relative error
is to be bounded. 
The analysis proceeds in three  steps: a step-by-step analysis of
the instructions of the program; an asymptotic analysis of an upper
bound on the
relative error as the precision tends to infinity (or equivalently, as $u$ tends to zero); an actual bound
on the quadratic term in the relative error when using the asymptotic
bound found for its linear part. We now review these steps in more
detail.

\subsection{Step-by-step analysis}\label{sec:step-by-step}
This is the
step where knowledge of
the FP numbers from Sections~\ref{subsec:FParith} and~\ref{sec:intro-bounds} is used. It
consists of an inductive construction of a system of polynomial
equations and a system of inequalities.

First, a
system~$I_0$ of linear inequalities is initialized with the known
information
on the ranges of the input variables, plus the inequalities
$0<u\le u_{\max}$ for the unit round-off (see 
\cref{sec:intro-bounds}). The initial system of polynomial equations
is the empty set: $S_0=\emptyset$.

The analysis proceeds instruction by instruction, constructing new
systems of equations and inequalities~$S_i$ and~$I_i$ from~$S_{i-1}$
and~$I_{i-1}$. Let~$t_i$ be the
right-hand side of the assignment in the $i$th instruction of the
algorithm from \cref{eq:instruction}. This step aims at bounding
the
rounding error in the assignment $v_i:=t_i$ by analyzing~$t_i$ using
the systems~$S_{i-1}$ and~$I_{i-1}$.

\paragraph{No error}
The best situation is when that error can be
determined to be~0. This happens either when $t_i$~is a constant
floating-point number that
does not depend on the input variables, when a multiplication by a power of $2$ is performed,  or when Sterbenz' 
\cref{lemma-sterbenz} or \cref{lemma-boldodaumasSQRT} are found to apply\footnote{That part of the
analysis is not fully implemented currently, but the user can give
this
information to our program.}. As we construct only polynomial
equations, the system~$S_i$ is then obtained by adding one of the
following
to~$S_{i-1}$:
\begin{equation}\label{eq-exact}
v_i=w_j\times w_k+w_\ell\quad\text{or}\quad v_iw_k=w_j\quad\text{or}\quad
v_i^2=w_j.
\end{equation}
In the last case, we also add~$v_i\ge0$ to~$I_{i-1}$ to obtain~$I_i$.
Otherwise, $I_i=I_{i-1}$. 
\paragraph{Absolute error}
The next best situation is when one can determine $\operatorname{ulp}
(t_i)$ exactly. For this, one determines the
minimal and maximal values of~$t_i$ given the equations and
inequalities in~$S_{i-1}$ and~$I_{i-1}$. If these bounds
allow to determine
$\ell_i=\lfloor\log_2|t_i|\rfloor$
exactly, then the system is enriched with an \emph{absolute error
bound}: in \cref{eq-exact}, $v_i$ is replaced by 
\begin{equation}\label{eq:absolute}
v_i+2^{\ell_i+1}\epsilon_iu,
\end{equation}
and the system of inequalities is complemented with 
$-1\le\epsilon_i\le 1$.

\paragraph{Relative error}
In the remaining cases, the systems are enriched with a relative
error bound:
\begin{align}
&&v_i&=(w_j\times w_k+w_\ell)(1+\epsilon_iu),& -b_{DK}&\le
\epsilon_i\le b_{DK}\label{eq:fma}\\
&\text{or}&
v_iw_k&=w_j(1+\epsilon_iu),& -b_{\div}&\le \epsilon_i\le b_
{\div},\label{eq:div}\\
&\text{or}&
v_i^2&={w_j}(1+\epsilon_iu)^2,&-b_{\sqrt{}}&\le \epsilon_i\le b_
{\sqrt{}},\quad v_i\ge0,\label{eq:sqrt}
\end{align}
and the relevant equation and inequality among
\begin{equation}\label{eq:bound-vars}
\begin{split}(1+u)b_{DK}=1,\quad b_\div=u-2u^2,\qquad\\
u(1+2u)b_{
\sqrt{}}^2-2
(1+2u)b_{\sqrt{}}+2=0,\quad b_{\sqrt{}}\ge 0,
\end{split}
\end{equation}
that come from
\cref{lemma:jeannerod-rump,lemma-knuth-dekker}.

\subsection{Polynomial optimization}
The end result of the step-by-step analysis is a system of polynomial
equations and a system of polynomial inequalities, where each
statement of the
algorithm brings one equation in one or two new variables (two being
the general case when the error is not zero) and a certain number of
inequalities.

At each stage of the step-by-step analysis and when bounding the
relative error on the last variable, we need to find the minimum or
maximum of a quantity defined by the system. Without loss of
generality and up to adding equations to the system, we can assume
that the quantity to be maximized or minimized is the last
variable~$v_m$ introduced in the system.

As mentioned in the introduction, this problem is
well-understood in terms of algorithms and complexity 
but also very expensive. In practice, even with efficient software such as Mohab Safey el Din's
\texttt{Raglib} library\footnote{\url{https://www-polsys.lip6.fr/~safey/RAGLib/}}, which relies on Faug\`ere's 
\texttt{FGb} library for
Gr\"obner bases\footnote{
\url{https://www-polsys.lip6.fr/~jcf/FGb/index.html}}, we could not
obtain our bounds directly using general-purpose approaches, except
for the simplest algorithms. 

Instead, we take an approach that exploits two characteristics
of our systems of polynomial equalities and inequalities:
\begin{enumerate}
    \item As the unit round-off $u$ is typically small compared to the other
variables of the program, the set of inequalities describes a
\emph{small tube-like
semi-algebraic set} around the exact value the variables take at~$u=0$;
    \item Being produced by the step-by-step analysis described above,
the set of polynomial equalities we start with has a \emph{triangular
structure}.
\end{enumerate}

\paragraph{Free and dependent variables.}
We distinguish
two types of variables in these systems: the \emph{free} variables are
the variables~$\epsilon_i$ encoding the
absolute or relative errors, the input variables and the unit
round-off~$u$; the \emph{dependent} variables are
those whose value is fixed once the free ones are fixed. Initially,
they
correspond to the variables~$v_i$ defined by the algorithm itself, as
well as the bound variables from \cref{eq:bound-vars}.

The optimization is a recursive process where at each step one
equation
is added to the system, making one of the free variables dependent
and maintaining the triangular nature of the system. At
the end, no free variables are left and the optimum is found among
finitely
many values. We now describe this process in more detail.

\paragraph{Gradient.}
As the system~$S_m$ obtained by the step-by-step
analysis is triangular, by repeated
use of
the chain rule, the gradient of the last variable~$u_m$ with respect to
the free variables is obtained as a vector of rational functions in
all variables. Equivalently, this can be computed by automatic
differentiation.

\paragraph{Recursive optimization.}
The optimization is a recursive search of extrema. It is
performed by looping
through the free variables and trying to detect whether
the
sign of the partial derivative with respect to one of
them can be
decided (see \cref{sec:sign-decision}). If this is the case, a new
equation is added to the system, setting this variable to
its extremal
value from \cref{eq:bound-vars} and one gets an expression in one variable less to be
optimized.
When the decision cannot be reached for any of the free variables,
then
we pick one of the variables~$\epsilon_i$ (heuristic choices are made)
and
use a branch-and-bound method, optimizing the polynomial in each of
the three cases $\partial u_m/\partial\epsilon_i=0$ or $\epsilon_i$
equal to its extremal values. The new systems have one variable less.

\subsection{Sign decisions}\label{sec:sign-decision}

Being able to quickly decide the sign of a polynomial in the domain of
optimization gives a substantial speed-up in the computation by
removing many branches in the optimization tree.
This is where we use the fact that several of the variables live in
a small domain. Several techniques are used to help this decision.

\paragraph{Factorization} of multivariate polynomials over the
rationals is well understood in theory
\cite{KaltofenLecerf2011,GathenGerhard2013} and works well in
practice. It is not needed by the algorithms but allows to reduce
the degrees and possibly the number
of variables of the polynomials whose signs have to be decided.

\paragraph{Interval arithmetic} can be applied since for each
variable, be it free or dependent, we have bounds at our disposal that
have been computed
during the step-by-step analysis. In some cases, this is sufficient to
decide the sign.

\paragraph{Small number of variables} occur near the bottom of the
recursion tree. A classical way to decide that
a polynomial in one variable does not vanish inside a given interval
is to use Sturm sequences. These are available in all major computer
algebra systems. 

For a univariate \emph{algebraic} expression~$A(x)$, e.g., a
combination of
$+,\times,\div,\sqrt{}$ applied to constants and one variable~$x$, it
is
easy to construct by induction a nonzero bivariate polynomial~$Q(x,y)$
such that $Q(x,A(x))=0$. Then a sufficient condition for $A$ not to
vanish in a given interval is that the constant term~$Q(x,0)$ does not
vanish there. This reduces the problem to that of a 
\emph{polynomial} in one variable. While this is only a
sufficient condition, this method saves a significant amount of time
when it applies.

\subsection{Regular chains}\label{subsec:regular-chains}
We take advantage of the \emph{triangular shape} of the
system constructed by the step-by-step analysis of 
\cref{sec:step-by-step} using \emph{regular chains}. These were
studied initially by Lazard~\cite{Lazard1991} and Kalkbrener~%
\cite{Kalkbrener1993} and a
large
number of efficient algorithms have been developed since~%
\cite{AubryLazardMoreno-Maza1999,ChenMoreno-Maza2011,ChenMoreno-Maza2012,ChenMoreno-Maza2016,BoulierLemaireMazaPoteaux2021,AlvandiChenMoreno-Maza2013}.

We first recall the basic notions. 
The polynomials belong to $\mathbb K[x_1,\dots,x_n]$ for an
algebraically closed field
$\mathbb K$.  $\mathbb K$ can be the field~$\mathbb C$ of
complex numbers, but usually, the coefficients of the polynomials
are rational functions in the free variables with rational
coefficients,
while the variables $x_i$ are the dependent variables.
An order~$x_n>\dots>x_1$ is fixed on the variables. The 
\emph{leading}
variable of a polynomial is the largest variable on which it
depends.
A set~$T$ of polynomials is 
\emph{triangular} when the polynomials have pairwise distinct leading
variables.

The coefficient of highest degree of a polynomial with respect
to its leading variable is called the
\emph{initial} of the polynomial.
The product of the initials of a
triangular set~$T$ is
denoted~$h_T$. Three geometric notions are relevant: the zero set~$V
(T)\subset\mathbb K^n$ of zeros of~$T$; the \emph{quasi-component} $W
(T)=V(T)\setminus V(h_T)$ and its Zariski closure~$
\overline{W (T)}$. 

Regular chains are defined inductively. An empty triangular set~$T$ is
a regular
chain. Otherwise, if $T_{\max{}}$ is the element of a triangular
set~$T$ with the
largest leading variable, $T$ is a regular chain when 
$T':=T\setminus\{T_{
\max{}}\}$ is a regular chain and the
the initial~$h_{\max}$ of the polynomial~$T_{\max{}}$ is regular in
the sense
that 
$\overline{W(T')}=\overline{W(T')\setminus V(h_{\max})}$. 

The following example illustrates all these notions.
\begin{example}Consider a program that would perform the sequence of
assignments
\begin{equation}\label{eq:assignments}
t_2:=\sqrt2t_1,\quad t_3:=t_1^2/t_2.
\end{equation}
For the second assignment not to raise an error, it is necessary
that~$t_2\neq0$. The computed values are such that the point~$
(t_1,t_2,t_3)$ is part of the solution set of the polynomial system
\[T=\{t_2^2-2t_1^2,t_3t_2-t_1^2\}.\]
This is a regular chain for the order $t_3>t_2>t_1$. Its solution set
is the union of two curves: $V(T)=\mathcal C_1\cup\mathcal
C_2$,
with
\[\mathcal C_1=\{(\pm\sqrt2t,2t,t)\mid t\in\mathbb C\},\qquad\mathcal
C_2=
\{(0,0,t)\mid
t\in\mathbb C\}.\]
The product of initials in~$T$ is~$h_T=t_2$ and thus the
quasi-component is just the first curve minus a point: $W(T)=\mathcal
C_1\setminus\{(0,0,0)\}$.
All values computed by \cref{eq:assignments} belong to~$W(T)$.
Finally, the Zariski closure recovers the missing point and
$\overline{ W(T)}=\mathcal C_1$.
\end{example}

\subsection{Application to Error Analysis}
The systems of equations constructed during the step-by-step analysis
of \cref{sec:step-by-step} are actually regular chains. The field one
starts with is the field of rational functions~$\mathbb K_0:=\mathbb Q
(u,i_1,\dots,i_\ell)$ in $u$ 
the unit round-off and $i_1,\dots,i_\ell$ the input variables.

Starting from
$I_{i-1}$ a regular chain in~$\overline{\mathbb K_{i-1}}
[v_1,\dots,v_{i-1}]$ for the order $v_1<\dots<v_{i-1}$, the chain is
enriched by the relevant equation from \cref{eq-exact}, 
\cref{eq:absolute} or
\cref{eq:fma,eq:div,eq:sqrt}, while the field $\mathbb K_i$ is 
$\mathbb K_{i-1}$ when it is detected that no error occurs and
$\mathbb K_{i-1}(\epsilon_i)$ otherwise. The new variable~$v_i$  is
set to be larger than $v_{i-1}$ in the ordering; the main variable of
the new
polynomial is thus~$v_i$, so that the system is triangular; its
initial is~1, except in the case of a division where it is~$w_k$. In
that situation, removing $V(w_k)$ is expected, as it corresponds to
a division by~0. Thus, computing quasi-components is desired in our
setting.

The construction of the regular chains used during the optimization
is different. This occurs during the step-by-step analysis to
determine
the extremal values of the variables or during the computation of
relative error bounds. During such an optimization, new
equations involve free variables that become dependent and are
removed
from the current field of coefficients. These new dependent variables
are appended to the list of dependent variables in \emph{descending
order}: if, before
the optimization, the system defines the variables~$v_m>v_
{m-1}>\dots>v_1$, during the optimization the variables are given by
a regular chain~$T$ that also
defines~$f_1>\dots>f_{k-1}$ with the $f_i$ formerly free variables
and $v_1>f_1$.

The introduction of a new variable~$f_k$ is through a polynomial~$p
(f_k)$. The equation $p(f_k)=0$ indicates
that a derivative is~0, or
that one of the~$\epsilon_i$, one of the input
variables,
or the
unit round-off~$u$ is assigned one of its possible extremal values.
A new field~$\tilde{\mathbb K}$ is defined with one variable less so
that the previous
one
is~$\tilde{\mathbb K}(f_k)$. 
At
this stage, the current regular chain~$T$ is 1-dimensional over the
current field,
as all dependent variables except~$f_k$ are given by a polynomial in
the
chain.
An important issue is that it is not sufficient to consider the
intersection~$W(T)\cap V(p(f_k))$, for which many algorithms are
available, but we need to know the
closure~$\overline{W(T)}\cap V(p(f_k))$. 
As an illustration of the difference, if the current
system is
\[T=[2t_3t_1^2-2(t_1+1)t_2-2t_1^2+3t_1+2,t_2^2-t_1-1]\]
with $t_3>t_2>t_1$, then the intersection~$W(T)\cap V(t_1)$ is empty,
whereas $\overline{W(T)}\cap V(t_1)=\{(0,1,11/8)\}$. The fact that our
construction always involves 1-dimensional regular chains lets us
compute these closures thanks to an algorithm by Alvandi et al.~%
\cite{AlvandiChenMoreno-Maza2013}, implemented in the function
\texttt{LimitPoints} of the Maple package
\texttt{RegularChains[AlgebraicGeometryTools]}.

\section{A gallery of algorithms for the hypotenuse}
\label{sec:algos-hypot}
Our running example is the calculation of the hypotenuse function  $
(x,y) \mapsto \sqrt{x^2+y^2}$.
We now present various algorithms for that function that we
use to illustrate our approach,  each with its own
motivation and merits.  The naive algorithm (Algorithm~\ref{NaiveHypot}) is simple, fast, and fairly accurate: when no
underflow/overflow occurs, the absolute error is bounded by $1.222~ \ulp$, as shown by Ziv~\cite{Ziv1999}, and the relative error is bounded by $2u$.
 However, Algorithm~\ref{NaiveHypot}  suffers from a serious
 drawback: intermediate calculations can underflow or overflow, even if
 $\sqrt{x^2+y^2}$ is far from the underflow or overflow thresholds.
 Such \emph{spurious} underflows or overflows can result in infinite
 or very inaccurate output. For instance, assuming we use the binary64
 (a.k.a. ``double precision'') format of the IEEE 754 Standard,
\begin{itemize}
    \item[] if $x = 2^{600}$ and $y = 0$ the returned result is $+\infty$ (because the computation of $x^2$ overflows) whereas the exact result is $2^{600}$;
    \item[] if $x = 65 \times 2^{-542}$ and $y = 72 \times 2^{-542}$ the returned result is $96 \times 2^{-542}$ whereas the exact result is $97 \times 2^{-542}$. 
\end{itemize}
The usual solution to overcome this problem is to \emph{scale} the input data,
i.e., to multiply or divide $x$ and $y$ by a common factor such that
overflow becomes impossible and underflow becomes either
impossible or harmless. This is for instance what Algorithm~\ref{simplest-scaling} below does, in a rather straightforward way.

To counter spurious overflow, one can divide both operands by the one
with the largest magnitude. This gives the well-known Algorithm~\ref{simplest-scaling}, which is used for example in Julia 1.1. Spurious underflow is not completely avoided, but if one of the input variables (say, $y$) is so small in magnitude before the other one that $|y/x|$ is below the underflow threshold, then $\sqrt{x^2+y^2}$ is approximated by $|x|$ with very good accuracy, so that the result returned by the algorithm is excellent.
Unfortunately, as we are going to see, Algorithm~\ref{simplest-scaling} is significantly less accurate than Algorithm~\ref{NaiveHypot}.
\begin{algorithm}[H]
 \caption{The simplest scaling for the \texttt{Hypot} function: divide
 both operands by the one with the largest magnitude. }
\label{simplest-scaling}
  \begin{algorithmic}[1] 
    \IF{$|x| < |y|$} \STATE{swap$(x,y)$} \ENDIF
    \STATE{$r \gets \RN(y/x)$}
    \STATE{$t \gets \RN(1+r^2)$}
    \STATE{$s \gets \RN(\sqrt{t})$}
    \STATE{$\rho_2 = \RN(|x| \cdot s)$}
  \end{algorithmic}
\end{algorithm}

Beebe~\cite{Beebe:2017} suggests several solutions to compensate for the loss of accuracy in Algorithm~\ref{simplest-scaling} (using an adaptation of the Newton-Raphson iteration for the square root). One of them is Algorithm~\ref{Beebe-alg} below, whose first $6$ lines are exactly the same as those of  Algorithm~\ref{simplest-scaling}. As shown in Section~\ref{sec:analysis-beebe}, it does more than just compensate for the loss of accuracy due to the scaling: it has a better final error bound than Algorithm~\ref{NaiveHypot}.

\begin{algorithm}[H]
 \caption{Improvement of Algorithm~\ref{simplest-scaling} suggested by Beebe~\cite{Beebe:2017}.}
\label{Beebe-alg}
  \begin{algorithmic}[1] 
    \IF {$|x| < |y|$}
    \STATE{swap$(x,y)$}
    \ENDIF
        \STATE{$r \gets \RN(y/x)$}
    \STATE{$t \gets \RN(1+r^2)$}
    \STATE{$s \gets \RN(\sqrt{t})$}
    \STATE{$\epsilon \gets \RN(t-s^2)$}
    \STATE{$c \gets \RN(\epsilon/(2s))$}
    \STATE{$\nu \gets \RN(|x| \cdot c)$}
    \STATE{$\rho_3 \gets \RN(|x| \cdot s + \nu)$}
  \end{algorithmic}
\end{algorithm}

Borges~\cite{Borges2020} presents several algorithms, including the very
accurate Algorithm~\ref{Borges-fused}
below. It does not address the spurious under/overflow problem (one
has to assume that a preliminary scaling by a power of $2$ has been done). The Fast2Sum and Fast2Mult algorithms called by Algorithm~\ref{Borges-fused} are well-known building blocks of computer arithmetic (see for instance~\cite{BoldoEtAl2023}). Here we just need to know that Fast2Mult$(x,y)$ (resp. Fast2Sum$(x,y)$) delivers a pair $(a,b)$ of FP numbers such that $a = \RN(xy)$ (resp. $a = \RN(x+y)$) and $b = xy-a$ (resp. $b = (x+y)-a$).
As shown in Section~\ref{sec:analysis:Borges}, this algorithm is
almost optimal in terms of relative error (we obtain a relative error
bound slightly above $u$).

\begin{algorithm}[H]
 \caption{Borges' corrected ``fused'' algorithm~\cite[Algorithm~5]{Borges2020}.}
\label{Borges-fused}
  \begin{algorithmic}[1] 
    \IF {$|x| < |y|$}
    \STATE{swap$(x,y)$}
    \ENDIF
    \STATE{$(s_x^h,s_x^\ell) \gets$ Fast2Mult$(x,x)$}
    \STATE{$(s_y^h,s_y^\ell) \gets$ Fast2Mult$(y,y)$}
    \STATE{$(\sigma_h,\sigma_\ell) \gets$ Fast2Sum$(s_x^h,s_y^h)$}
    \STATE{$s \gets \RN\left(\sqrt{\sigma_h}\right)$}
    \STATE{$\delta_s \gets \RN(\sigma_h - s^2)$}
    \STATE{$\tau_1 \gets \RN(s_x^\ell + s_y^\ell)$}
    \STATE{$\tau_2 \gets \RN(\delta_s + \sigma_\ell)$}
    \STATE{$\tau \gets \RN(\tau_1 + \tau_2)$}
    \STATE{$c \gets \RN(\tau/s)$}
    \STATE{$\rho_4 \gets \RN(c/2 + s)$}
\end{algorithmic}
\end{algorithm}

Kahan~\cite{Kahan1987} gives the fairly accurate  Algorithm~\ref{alg-Kahan} below. It avoids spurious underflows and overflows. Although we do not discuss these matters here, it also has the advantage of correctly setting the various IEEE 754 exception flags (underflow, overflow, inexact).  It has the kind of size where automation of the analysis becomes a necessity if one is not satisfied with loose error
bounds. 
\begin{algorithm}[t]
 \caption{Kahan's \texttt{hypot} Algorithm~\cite[Algorithm~\texttt{CABS}]{Kahan1987}.
 In this presentation, it requires the availability of an FMA instruction. We assume $0 \leq y \leq x$. The algorithm uses the precomputed constants $R_2 = \RN(\sqrt{2})$, $P_h = \RN(1+\sqrt{2})$, and $P_\ell =  \RN(1+\sqrt{2}-P_h)$.}
\label{alg-Kahan}
 \begin{algorithmic}[1]
 \STATE{$\delta \gets \RN(x-y)$}
\IF{$\delta > y$}
\STATE{$r \gets \RN(x/y)$}
\STATE{$t \gets \RN(1+r^2)$}
\STATE{$s \gets \RN(\sqrt{t})$}
\STATE{$z \gets \RN(r+s)$}
\ELSE
\STATE{$r_2 \gets \RN(\delta/y)$}
\STATE{$tr_2 \gets \RN(2r_2)$}
\STATE{$r_3 \gets \RN(tr2 + r_2^2)$}
\STATE{$r_4 \gets \RN(2+r_3)$}
\STATE{$s_2 \gets \RN(\sqrt{r_4})$}
\STATE{$d = \RN(R_2 + s_2)$}
\STATE{$q = \RN(r_3/d)$}
\STATE{$r_5 \gets \RN(P_\ell+q)$}
\STATE{$r_6 \gets \RN(r_5+r_2)$}
\STATE{$z \gets \RN(P_h+r_6)$}
\ENDIF
\STATE{$z_2 \gets \RN(y/z)$}
\STATE{$\rho_5 \gets \RN(x+z_2)$}
 \end{algorithmic}
\end{algorithm}

These algorithms are sorted by order of complexity of their analysis.
We now deal with them in order, introducing the new problems and
solutions one at a time.

\section{Algorithm~\ref{NaiveHypot}: Straightforward Analysis}\label{sec:warm-up}

As mentioned in the introduction,  the
relative error of Algorithm~\ref{NaiveHypot} is bounded by $2u$ and that
bound is asymptotically optimal.
We now show how this result is derived by our approach, and how a small
term $\frac{5}{4}u^2$ can be subtracted from the previous bound.

\paragraph{Step-by-step analysis}
The result of a step-by-step analysis was given in \cref{eq-algo1}. It follows that
the relative error $\rho/\sqrt{x^2+y^2}-1$ equals
\[R=\frac{\sqrt{\left(x^2(1+u\epsilon_{s_x})+y^2(1+u\epsilon_
{s_y})\right)
(1+u\epsilon_\sigma)}\,(1+u\epsilon_\rho)}{\sqrt{x^2+y^2}}-1,\]
where each of $\epsilon_{s_x},\epsilon_{s_y},\epsilon_\sigma$ has
absolute value bounded by~$1/
(1+u)$, by \cref{rel-err-v}, while $\epsilon_\rho$ is bounded by the
Jeannerod-Rump bound of
\cref{err-sqrt-jr}. 
The next step of the analysis is to find the maximal value of~$|R|$
when all the
$\epsilon$ variables range in these intervals. Note that all these bounds are smaller than~1.

\paragraph{Upper bound}
Since $0<u<1$, the expression above is an
increasing function
of each of these~$\epsilon_i$ in their intervals. Therefore its
minimum is reached when they are all simultaneously equal
to their minimum, and similarly for the maximal value. When 
$\epsilon_{s_x}=\epsilon_{s_y}=\epsilon_\sigma$ with $\epsilon$ their
common value, the
expression of~$R$ simplifies to
\[(1+u\epsilon)(1+u\epsilon_\rho)-1=u
(\epsilon+\epsilon_\rho)+u^2\epsilon\epsilon_\rho.\]
This shows that the \emph{absolute value} of the relative error~$R$ is
maximal when $\epsilon$
and $\epsilon_\rho$ both reach their maximal value, giving
\[
|R|\le S(u):=\frac{1+2u}
{1+u}\left(2-\frac1{\sqrt{1+2u}}\right)-1= 
\frac{1+3u-\sqrt{1+2u}}{1+u}.
\]
As $u\rightarrow0$, 
\[S(u)  = 2u - \frac{3}{2}u^2 + u^3 + O(u^4),\]
showing the linear term in the bound.
If only the linear term $2u$ is needed, using the simple
bound from
\cref{rel-err-u} is sufficient. The more refined estimates are only required to
control the difference between that linear term and the actual error.

\paragraph{Optimal quadratic term}
Next, consider 
\[y(u)=\frac{S(u)-2u}{u^2},\]
whose maximum for $u\in(0,u_{\max}]$ gives an upper bound on the
quadratic term of
the error bound. This function is $C^\infty$ for $u>-1/2$ and the Taylor
expansion above shows that $y(0)=-3/2$ and $y'(0)=1$.
From the explicit expression of~$S$, it is easy to see that 
$y'>0$ in the interval $[0,1/4]$, implying that for $u_{\max}\le 1/4$
its
maximum~---~the bound we are after~---~is reached at~$u_{\max}$. Thus
we have proved the following.
\begin{theorem}
\label{th:naive-alg}
Barring underflow and overflow, the relative error~$R$ of 
Algorithm~\ref{NaiveHypot} for $u\le 1/4$ satisfies
\[|R|\le \frac{1+3u-
\sqrt{1+2u}}{1+u}=2u+\kappa u^2,\quad \kappa<-5/4.\]
\end{theorem}
The bound given in the introduction and obtained by our
program for $u=1/4$ can be seen to be $y(1/4)$.
The bound $y(1/4)<-5/4$
follows by numerical approximation. Other values of the quadratic
bound are given in \cref{table:quad-algo1}.

\begin{table}
\begin{center}
\renewcommand{\arraystretch}{1.2}
\begin{tabular}{cl}
\toprule
lower bound on $p$ & upper bound on the error\\
\midrule
8&$2u-\frac32u^2+4\cdot10^{-3} u^2$\\
11&$2u-\frac32u^2+5\cdot 10^{-4} u^2$\\
24&$2u-\frac32u^2+6\cdot 10^{-8} u^2$\\
53&$2u-\frac32u^2+2\cdot 10^{-16} u^2$\\
113&$2u-\frac32u^2+10^{-34} u^2$\\
\bottomrule
\end{tabular}
\renewcommand{\arraystretch}{1.0}
\end{center}
\caption{Quadratic error bound for Algorithm~\ref{NaiveHypot}}
\label{table:quad-algo1} 
\end{table}

\paragraph{Proof based on polynomials}
In preparation for more complicated examples where the analysis of
the function~$S(y)$ is not as straightforward as here, we show
how
computer
algebra could be used to show that $y'>0$ in the interval~$[0,1/4]$.
The starting point is to perform a simple manipulation (automated in
the Maple \texttt{gfun} package) showing that~$y$
satisfies the quadratic equation
\[E(y,u)=u^2(1+u)^2y^2-2(1-u^2)(1+2u)y+(1+2u)(2u-3)=0.\]
By continuity, $y'$ is positive in a neighborhood of~$u=0$. It can only change
sign at a point~$u$ where $y'=0$, but then 
\[\frac{\partial E}{\partial u}+\frac{\partial E}{\partial y}y'=0\]
implies $\partial
E/\partial u(y,u)=0$.
Eliminating $y$ between $E=0$ and $\partial E/\partial u=0$ by means of a
resultant shows that in turn, this implies
\[(1+u)(4+21u+32u^2+12u^3-8u^4)=0.\]
Now, a simple interval analysis shows that this
polynomial does not vanish for $u\in[0,1/2]$. Thus $y'$ does not vanish
in this interval and $y'(0)>0$ shows that $y$ is an increasing function.

Then, a direct proof of the bound
$-5/4$ consists in using
continuity again, considering the univariate polynomial $E
(-5/4,u)$ and observing its absence of root in $[0,1/2]$, which shows
that $-5/4$ is never reached by~$y(u)$ in this interval.

\section{Algorithm~\ref{simplest-scaling}: Exploit Absolute Errors}
\label{sec:analysis-algo2}
A direct analysis using only the simple relative error bound 
\cref{rel-err-u}
leads to a bound on the relative error with a linear term of order $3u$. The
refined estimates on the relative error from \cref{sec:intro-bounds} do not
improve that linear term. A step-by-step analysis of the program
yields a better estimate by
studying more carefully the ranges of the intermediate variables.

We assume $0 \leq y \leq x$ (i.e., we consider that if needed, the swap of Line 2 of the algorithm has already been done). 

\subsection{Automatic Analysis}
The result of our automation is as follows.
\begin{Verbatim}
> Algo2:=[Input(x=0..2^16,alpha=0..1),y=RN(alpha*x,0),
> r=RN(y/x),t=RN(1+r^2),s=RN(sqrt(t)),rho=RN(x*s)]:
> BoundRoundingError(Algo2);
\end{Verbatim}
\[{\color{blue}\frac52\_u+\frac38\_u^2}\]
The assumption $0\le y\le x$ is encoded by introducing a variable $\alpha\in
[0,1]$ and stating that $y=\alpha x$ without an error, by using a second
argument to \texttt{RN} that encodes special knowledge on the error
for a specific operation.

\subsection{Result}
The output of our program is an indication pointing to the following result,
proved below.
\begin{theorem}
\label{theorem-tighter-alg2}
The relative error of Algorithm~\ref{simplest-scaling} is less than or equal to
\[
R_5(u) = \frac{(1+2u)\sqrt{1+u}-1+2u^2}{1+u} \leq  \frac{5}{2} u + \frac{3}{8} u^2.
\]
\end{theorem}
\begin{proposition}
The bound of Theorem~\ref{theorem-tighter-alg2} is asymptotically optimal.
\end{proposition}
\begin{proof}
Choose $x = 2^p-1$ and
\[
y = \RN\left( (2^p-1) \cdot \sqrt{7} \cdot 2^{-p/2} \cdot (1-u)^2 \right).
\]
As soon as $u \leq 1/32$, following the steps of the algorithm
shows that its output is $\rho_2 = 2^p$. It follows 
that $\sqrt{x^2+y^2}-\rho_2 \to \frac{5}{2}$ as $u \to 0$, so that the relative error is asymptotically equivalent to $5u/2$.
\end{proof}
 
\noindent The following example illustrates the sharpness of the bound:
   in the binary64 format, with 
   \[x = 9007199254740991\quad\text{and}\quad y =
  {8425463406411589} \times 2^{-25},\] the relative error is
   $2.49999999999999558648 u$.

We now detail the steps of the analysis leading to the proof of 
\cref{theorem-tighter-alg2}, following the approach
outlined in \cref{app-ca}.

\subsection{Step-by-step analysis}
This stage translates the program into a polynomial system with bounds on
intermediate error-variables.
The first two steps give the equalities
\[y=\alpha x,\quad rx=y(1+u\epsilon_r),\]
with $|\epsilon_r|\le 1-2u$ by the Jeannerod-Rump bound \cref{err-div-jr}.

Next, since $r$ is the rounded value of
$\alpha\in[0,1]$ and the function $\RN$ is increasing, $r \leq \RN(1) =
1$. This implies that $1\le 1+r^2 \leq 2$ and therefore also $1\le\RN(1+r^2)\le
2$. Hence the rounding \emph{absolute} error committed at Line~5 
of Algorithm~\ref{simplest-scaling} is less than $u$ (half the  distance between two consecutive FP numbers between $1$ and $2$).
The same reasoning applies to~$\sqrt{t}$. This gives
\[t=1+r^2+u\epsilon_t,\quad s=\sqrt{t}+u\epsilon_s,\]
with $|\epsilon_s|$ and $|\epsilon_t|$ bounded by~$1$.

The last step is analyzed as before, giving
\[\rho=xs(1+u\epsilon_\rho),\]
with $|\epsilon_\rho|\le 1/(1+u)$, by \cref{rel-err-v}.

It follows from these equations that the relative error $|\rho/
\sqrt{x^2+y^2}-1|$ is upper bounded by $|R|$, where
\[R=\frac{\left(\sqrt{1+\alpha^2
(1+u\epsilon_r)^2+u\epsilon_t}+u\epsilon_s\right)(1+u\epsilon_\rho)}{\sqrt{1+\alpha^2}}-1.\]

\subsection{Upper bound}
As in the analysis of \cref{sec:warm-up}, for $u\le 1/4$, since all the
$|\epsilon_i|$ are bounded by~1, $R$ above is an increasing function
of all the~$\epsilon_i$ in their intervals. 
Writing
\[\epsilon_r=\epsilon(1-2u),\quad
\epsilon_s=\epsilon,\quad
\epsilon_t=\epsilon,\quad
\epsilon_\rho=\frac\epsilon{1+u},
\]
it follows that $R$ is bounded between the values taken by
\[\tilde R(\alpha,\epsilon):=
\frac{\left(\sqrt{1+\alpha^2
(1+\epsilon u(1-2u))^2+u\epsilon}+u\epsilon\right)
\left(1+\epsilon\frac u{1+u}\right)}{\sqrt{1+\alpha^2}}-1
\]
at $\epsilon=-1$ and at $\epsilon=1$. The variation of these bounds
with respect to~$\alpha$ is dictated by
\begin{multline*}\frac{\partial\tilde R}{\partial\alpha}=
\left(1+\epsilon\frac u{1+u}\right)
\frac\alpha{(1+\alpha^2)^{3/2}}\\
\times\frac{(1+\epsilon u(1-2u))^2-1-u\epsilon-u\epsilon
\sqrt{1+\alpha^2
(1+\epsilon u(1-2u))^2+u\epsilon}}{\sqrt{1+\alpha^2
(1+\epsilon u(1-2u))^2+u\epsilon}}.
\end{multline*}
With $u\le 1/4$, $\epsilon\in\{-1,1\}$ and $\alpha\in[0,1]$,
the signs of each of these factors can be analyzed directly. The
first two are nonnegative, as is the denominator of the third one. 
Its numerator rewrites as
\[u\epsilon\left(1-(4-\epsilon)u-4\epsilon u^2(1-u)-\sqrt{1+\alpha^2
(1+\epsilon u(1-2u))^2+u\epsilon}\right).\]
The first factor is positive, the second one has the sign of~$\epsilon$ and the
last one is negative. Indeed, the argument of the square root is increasing with
$\epsilon$ and $\alpha$, so that its minimum is larger than $\sqrt{1-u}$, itself
larger than $1-u$ for $u\in[0,1]$. Thus this last factor is upper bounded by
\begin{multline*}
-1+u+1-(4-\epsilon)u-4\epsilon u^2(1-u)\le -3u+\epsilon u-4\epsilon u^2
(1-u)\\ 
\le -2u+4u^2(1-u)\le -u<0.
\end{multline*}
In conclusion, at $\epsilon=-1$, $\tilde R$ is negative and increasing with
$\alpha$, while at
$\epsilon=1$, it is positive and decreasing with $\alpha$. Thus in both cases,
its absolute value is maximal at~$\alpha=0$. Therefore~$R$ is bounded between
the values taken by
\[\tilde R(0,\epsilon):=
{\left(\sqrt{1+u\epsilon}+u\epsilon\right)
\left(1+\epsilon\frac u{1+u}\right)}-1
\]
at $\epsilon=-1$ and at $\epsilon=1$. It can now be observed that this is
maximal in absolute value at~$\epsilon=1$ so that we have obtained the upper
bound of the theorem
\[|R|\le S(u):=
\frac{(1+2u)\sqrt{1+u}-1+2u^2}{1+u}.\]
As $u\rightarrow0,$
\[S(u)=\frac52 u+\frac38u^2-\frac9{16}u^3+O(u^4).\]

\subsection{Optimal quadratic term}
Set $y(u)=({S(u)-\frac52u})/{u^2}$,
which satisfies the quadratic equation
\[E(y,u)=4u^2(1+u)^2y^2+4(1+u)(u^2+5u+2)y+u^2-6u-3=0.\]
We now show that $y'<0$ in the interval $[0,u_{\max}]$, implying that
the
maximum of~$y$ is reached at $u=0$, where it is $3/8$, as can be seen
from the Taylor
expansion above.

By continuity, $y'$ is negative in a neighborhood of~$u=0$. The vanishing
of~$y'$ can only occur at a zero of the resultant of~$E$ and $\partial
E/\partial u$, which implies
\[(1+u)^3(1+2u)^2(u^3-16u^2-22u-9)=0.\]
A simple analysis (e.g., by Sturm sequences) shows that this does not vanish
for $u$ in $[0,1]$. This proves that $y'<0$ in that region and concludes the
proof of the theorem.

 \section{Algorithm~\ref{Beebe-alg}: Split Argument Interval}
 \label{sec:analysis-beebe}

 In the analysis of Algorithm~\ref{Beebe-alg}, a new technique is necessary in
 order to obtain an optimal linear bound: splitting the domain of the
 parameters into two subdomains. The location of the split
 follows from the analysis and is presented as a suggestion by the
 implementation.
  
As in the previous analysis, we assume that, if needed, the swap of
Line~2 of the algorithm already took place, and we assume without loss
of generality that the operands are positive,  so that $0 \leq y \leq
x$.  We also assume $p \geq 4$ (i.e., $\umax =  1/16$). Also, as in
the previous algorithm, we set $\alpha = {y}/{x}$ and discuss
according to its value.

\subsection{Automatic analyses}
A first analysis is misleading:
\begin{Verbatim}
> Algo3:=[Input(x=0..2^16,alpha=0..1),y=RN(alpha*x,0),
> r=RN(y/x),t=RN(1+r^2),s=RN(sqrt(t)),
> epsilon=RN(t-s^2,0),c = RN(epsilon/2/s),nu = RN(x*c),
> rho=RN(nu+x*s)]:
> BoundRoundingError(Algo3);
\end{Verbatim}
\[\color{blue}
\frac{7}{4} \textit{\_u} +\left(\sqrt{2}-\frac{33}{32}\right) \textit{\_u}^{2}
\]
Recall that the second argument~0 of the rounding function \texttt{RN}
indicates that this operation is
exact. This is used as in Algorithm~\ref{simplest-scaling} to indicate
$y\le x$ and here also for the computation of~$\epsilon$, in view of
\cref{lemma-boldodaumasSQRT}.

While the result above is correct, it is pessimistic. 
A more careful analysis detailed below shows that it is beneficial to
analyze  the cases $\alpha\le1/2$ and $\alpha\ge 1/2$ separately. This
is hinted at by the implementation at a sufficiently high verbosity
level: it outputs
{\color{blue}
\begin{Verbatim}
getAbsoluteError: Splitting at r = 1/2 may help improve bounds
\end{Verbatim}
}
Once
this information is fed into our code (by changing the input range
\verb|0..1| into \verb|0..1/2| or \verb|1/2..1|, leading to the
variants \texttt{Algo3\_1} and \texttt{Algo3\_2} of the algorithm),
more precise
estimates are obtained:
\begin{Verbatim}
> BoundRoundingError(Algo3_1);BoundRoundingError(Algo3_2);
\end{Verbatim}
\begin{gather*}
\color{blue}
\frac{8 \textit{\_u}}{5}+\left(-\frac{20992}{5}+\frac{1447357555743 
\sqrt{5}}{1804709172500}+\frac{183265505148 \sqrt{21361601}\, \sqrt{5}}{451177293125}\right) \textit{\_u}^{2}\\
\color{blue}\frac{8 \textit{\_u}}{5}+\left(-\frac{20992}{5}+
\frac{22912648671 \sqrt{4173}}{352550900}+\frac{1413709473 \sqrt{5}}{1762754500}\right) \textit{\_u}^{2}
\end{gather*}
The linear term is improved from~$7/4=1.75$ down to~$8/5=1.6$.
The coefficients of~$\_u^2$ are approximately~$1.325$ and~$1.734$, so 
that the second one dominates. A finer analysis of the rounding error
on~$c$, detailed in the proof of \cref{th-sharp-alg-beebe} below,
shows that its absolute error is bounded by $u^2/2$. This information
can be incorporated into the statement of the algorithm and
leads to a further improvement of the error bound:
\begin{Verbatim}
> Algo3_2_better:=[Input(x=0..2^16,alpha=1/2..1),
> y = RN(alpha*x,0), r = RN(y/x), t = RN(1+r^2), s = RN(sqrt(t)),
> epsilon = RN(t-s^2,0), c = RN(epsilon/2/s,'absolute'(_u^2/2)),
> nu = RN(x*c), rho=RN(nu+x*s)];
> BoundRoundingError(Algo3_2_better);
\end{Verbatim}
\[\color{blue}\frac{8 \textit{\_u}}{5}+\left(-\frac{20992}{5}+\frac{5728153728 \sqrt{4173}}{88137725}+\frac{4126452 \sqrt{5}}{6779825}\right) \textit{\_u}^{2}\]
The coefficient of~$\_u^2$ is approximately $1.295$. We prove the
bound $1.4$ below, together with more precise ones for specific values
of~$p$ (see \cref{table:quad-algo3}).
\subsection{Result}
The analysis detailed in the following subsections follows the steps
of the automatic analysis, combined with a refined analysis of the
error on~$c$. It results in the following.
\begin{theorem}
\label{th-sharp-alg-beebe}
Assuming $u \leq 1/16$ (i.e., $p \geq 4$) and that an FMA instruction is
available, the relative error of Algorithm~\ref{Beebe-alg} is bounded by
\begin{small}
\begin{align*}
\chi_4(u)   &={\left(1+2 u \right) \sqrt{\frac{1+
u/5}{1+u}}}-1+u^{2}\frac{ \left(1+2 u\right)^{2}}{\left(1+u \right)^{2}}
\left(
\frac{\sqrt{5}}{5}+\frac{1}{\frac{{5}\, \sqrt{\left(1+u
\right) \left(1+\frac{u}{5}\right)}}{2}-u}
+\frac{2 \sqrt{5}}{5 \left(1+2 u \right)}\right),\\
&=\frac{8}{5} u +\left(\frac{3 \sqrt{5}}{5}-\frac{2}{25}\right) u^{2}+\left(\frac{116}{125}+\frac{14 \sqrt{5}}{25}\right) u^{3}+\mathrm{O}\! \left(u^{4}\right)\\ 
&\simeq1.6u+1.26u^2+O(u^3)\\ 
&\le \frac85u+\frac{7}{5}u^2,\qquad u\in[0,1/16].
\end{align*}
\end{small}
Moreover, the bound $\chi_4(u)$ is an increasing function of $u$.
\end{theorem}
\noindent Quadratic bounds for several values of~$p$ are given in 
\cref{table:quad-algo3}.

\begin{table}
\begin{center}\begin{tabular}{cl}
\toprule
lower bound on $p$ & upper bound on the error \\
\midrule
4&$1.6u+1.392u^2$\\
5&$1.6u+1.329u^2$\\
6&$1.6u+1.296u^2$\\
7&$1.6u+1.279u^2$\\
8&$1.6u+1.271u^2$\\
\bottomrule
\end{tabular}\end{center}
\caption{Quadratic error bound for Algorithm~\ref{Beebe-alg}} 
\label{table:quad-algo3} 
\end{table}

We are not able to prove that the linear term~$(8/5)u$ is
asymptotically optimal. Still, it is  sharp, in the
sense
given in \cref{sec:kindofbounds}, as shown by the following examples:
\begin{itemize}
   \item if $u = 2^{-53}$ (binary64 format), then error $1.5999739 u$ is attained with $x = 8056283928243985$ and $y = 4028141964171097$;
   \item if $u = 2^{-113}$, which corresponds to the binary128 (a.k.a. quad-precision) format of IEEE 754, then error $1.5999999648 u$, is attained with 
   \begin{align*}x &= 9288262988033986935972257666807793\\
   \intertext{and}
   y &=
   4644131494016993467987768200983857.
   \end{align*}
\end{itemize}

\subsection{Step-by-step analysis}
The core idea of the algorithm comes
from Newton's iteration, which translates as
\begin{equation}\label{eq:Newton}
\frac\epsilon{2s}+s=\frac{t-s^2}{2s}+s=\sqrt{t}+\frac{(s-
\sqrt{t})^2}{2s},
\end{equation}
where the last term exhibits the quadratic convergence.
The beginning of the analysis is the same as in 
\cref{sec:analysis-algo2}, with absolute error bounds for~$t$ and~$s$,
leading to the system
\begin{equation}\label{eq:algo3-part1}
y=\alpha x,\quad \underline{rx=y(1+u\epsilon_r)},\quad
t=1+r^2+u\epsilon_t,\quad s=\sqrt{t}+u\epsilon_s,\end{equation}
with $|\epsilon_r|\le 1-2u$, $|\epsilon_s|$ and $|\epsilon_t|$ bounded
by~1. In view of the upcoming discussion on~$\alpha$, the analysis in
this section does not make use of the underlined equation.
The next step gives $\epsilon=t-s^2$ (no rounding error by \cref{lemma-boldodaumasSQRT}). 

\paragraph{Detailed analysis for $c$}
From 
\cref{eq:Newton}, it follows that
\[\left|\frac\epsilon{2s}\right|=\left|-u\epsilon_s+\frac{u^2\epsilon_s^2}
{2s}\right|\le u+\frac{u^2}{2}.\]
If $\left|{\epsilon}/{(2s)}\right| \leq u$ then the error committed by rounding
$\frac{\epsilon}{2s}$ to nearest is less than $u^2/2$.
If $\left|{\epsilon}/{(2s)}\right| > u$, then, since the floating-point number
immediately above $u$ is $u + 2u^2$, the upper bound above implies that
$\RN\left({\epsilon}/{(2s)}\right) = \pm u$, so that again the rounding error is
less than $u^2/2$.
Therefore, in all cases, $|c| \leq u$ and
\begin{equation}\label{eq:algo3-part3}
c=\frac\epsilon{2s}+\epsilon_c\frac{u^2}2,
\end{equation}
with $|\epsilon_c|\le1$\footnote{This type of analysis is not automated yet.
By default, our program uses the coarser bound $c=\epsilon
(1+u\epsilon_c)/(2s)$, with
$|\epsilon_c|\le 1-u$. It can be given this extra information
through the two-arguments \texttt{RN}. The difference only
affects the quadratic term in the
final bound.}.
The next steps give 
\begin{equation}\nu=xc(1+u\epsilon_\nu),\quad \rho=(\nu+xs)(1+u\epsilon_\rho),
\label{eq:algo3-part2}\end{equation}
with $|\epsilon_\nu|$ and
$|\epsilon_\rho|$ bounded by $1/(1+u)$ by \cref{rel-err-v}. 
The result of this analysis is the following formula.
\begin{lemma}\label{lemma:rel-err-3}
The relative error of Algorithm~\ref{Beebe-alg} is bounded by
\begin{align}
R&=\sqrt{1+\frac{r^2-\alpha^2}
{1+\alpha^2}}
\sqrt{1+\frac{u\epsilon_t}{1+r^2}}\label{eq:defR-algo3}\\ 
&\qquad\qquad\times
\left(1+\frac{u^2}{2\sqrt{t}}
((\epsilon_c+\epsilon_s^2/s)(1+u\epsilon_\nu)-2\epsilon_s\epsilon_\nu)
\right)
(1+u\epsilon_\rho)-1,\notag\\
&=\frac{r^2-\alpha^2+u\epsilon_t}
{2(1+\alpha^2)}+u\epsilon_\rho+O
(u^2),\quad u\rightarrow0.\notag
\end{align}
Moreover, in this formula, $|\epsilon_s|,|\epsilon_t|,|\epsilon_c|$ are bounded
by~1 and $|\epsilon_\nu|$ and $|\epsilon_\rho|$ by $1/(1+u)$. 
\end{lemma}
\begin{proof}
This is obtained from the previous equations by a sequence of
rewriting
operations
\begin{align}
\rho&=(\nu+xs)(1+u\epsilon_\rho),\notag \\
&=x\left((-u\epsilon_s+\frac{u^2}2(\epsilon_c+\epsilon_s^2/s))
(1+u\epsilon_\nu)+\sqrt{t}+u\epsilon_s\right)
(1+u\epsilon_\rho),\notag \\
&=x\left(\sqrt{t}+\frac{u^2}2
((\epsilon_c+\epsilon_s^2/s)(1+u\epsilon_\nu)-2\epsilon_s\epsilon_\nu)
\right)
(1+u\epsilon_\rho),\label{eq:algo3-simple-rho} \\
&=x\sqrt{1+r^2}\sqrt{1+\frac{u\epsilon_t}{1+r^2}}
\left(1+\frac{u^2}{2\sqrt{t}}
((\epsilon_c+\epsilon_s^2/s)(1+u\epsilon_\nu)-2\epsilon_s\epsilon_\nu)
\right)
(1+u\epsilon_\rho)\notag,
\end{align}
whence the relative error of \cref{eq:defR-algo3}.
\end{proof}
The rest of the proof of \cref{th-sharp-alg-beebe} can be found in
\cref{app-proof-beebe}.

\section{Algorithm~\ref{Borges-fused}: Limited Human Proof}
\label{sec:analysis:Borges}
\subsection{Automatic Analyses}
This algorithm is beginning to be large for the current version of our
code, but can still be analyzed automatically:

\begin{Verbatim}
> Algo4:=[Input(x=0..2^16,alpha=0..1),y=RN(alpha*x,0),sxh=RN(x^2),
> sxl=RN(x^2-sxh,0),syh=RN(y^2),syl=RN(y^2-syh,0),
> sigmah=RN(sxh+syh),sigmal=RN(sxh+syh-sigmah,0),
> s=RN(sqrt(sigmah)),deltas=RN(sigmah-s^2,0),tau1=RN(sxl+syl),
> tau2=RN(deltas+sigmal),tau=RN(tau1+tau2),c=RN(tau/s),
> rho=RN(c/2+s)]:
> BoundRoundingError(Algo4);
\end{Verbatim}
\[\color{blue}
\textit{\_u} +\left(\frac{585981351743 \sqrt{66}}{1142440000}-4160\right) \textit{\_u}^{2}
\]
(Recall that our default value of $u_{\max{}}$ is~$1/64$.) 
The bound on the quadratic term is approximately~$6.989$.
Looking at the computation more closely, we make the following.
\begin{conjecture}\label{conj-beebe} For $u\le 1/4$, the relative
error of
\cref{Borges-fused} is bounded by
\begin{align*}
\phi(u)&:=\frac{(2+8 u +27 u^{2}+51 u^{3}+24 u^{4}-40 u^{5}-48 u^
{6}-16 u^
{7})\sqrt{1+2u}}{2(1+u)^4}-1\\
&\le u+7u^2.
\end{align*}
\end{conjecture}
\subsection{Result}
Our result on this algorithm is not as tight as \cref{conj-beebe},
but it does not exceed its value too much.
\begin{theorem}\label{thm-alg-Borges}
Barring overflow and underflow, as soon as
$p\ge 4$, the relative error of 
\cref{Borges-fused} is bounded by
\[u+(7+\kappa)u^2,\qquad\text{with}\quad 
\kappa\le\begin{cases}
21.4,&\quad u \leq 2^{-4},\\
6.1,&\quad u \leq 2^{-5},\\
2.5,&\quad u \leq 2^{-6},\\
1.2,&\quad u \leq 2^{-7},\\
0.6,&\quad u \leq 2^{-8},\\
7\,10^{-2},&\quad u \leq 2^{-11},\\
8\,10^{-6},&\quad u \leq 2^{-24},\\
2\,10^{-14},&\quad u \leq 2^{-53},\\
2\,10^{-32},&\quad u \leq 2^{-113}.
\end{cases}\]
\end{theorem}
\noindent The (tedious) proof is given in \cref{app-proof-thm-Borges}.
\begin{remark}
The bound given by Theorem~\ref{thm-alg-Borges} is asymptotically equivalent to~$u$. Such a bound is asymptotically optimal
(for any algorithm), as shown by the
following examples:
\begin{itemize}
    \item If $p$ is odd, then for $x=1$ and $y = 2^{(-p+1)/2} = \sqrt{2u}$, we have
  \[
   1+u-\frac{u^2}{2} < \sqrt{x^2+y^2} < 1+u-\frac{u^2}{2} + \frac{u^3}{2},
    \]
    so that $\sqrt{x^2+y^2}$ is at a distance at least $u - \frac{u^2}{2}$ from a FP number,  which implies that the relative error committed when evaluating it with any algorithm   is larger than
    \[
    \frac{u-\frac{u^2}{2}}{1+u-\frac{u^2}{2} + \frac{u^3}{2}} = u - \frac{3}{2}u^2+\mathcal{O}(u^3).
    \]
   \item If $p$ is even, then for $x=1$ and $y = \left\lceil \sqrt{2} \cdot 2^{p-1} \right\rceil \cdot 2^{-3p/2+1}$, we have
   \[
   1+u-\frac{u^2}{2}< \sqrt{x^2+y^2} < 1 + u + 2\sqrt{2}u^2+2u^3,
   \]
   so that $\sqrt{x^2+y^2}$ is at a distance at least $u-2\sqrt{2}u^2-2u^3$ from a FP number, which implies that the relative error committed when evaluating it with any algorithm  is larger than
   \[
   \frac{u-2\sqrt{2}u^2-2u^3}{1 + u + 2\sqrt{2}u^2+2u^3} = u - (1+2\sqrt{2})u^2 +\mathcal{O}(u^3).
   \]
   \end{itemize}
\end{remark}

\section{Algorithm~\ref{alg-Kahan}: Computer-aided analysis}
\label{sec:analysis-Kahan}
We first briefly explain the main ideas behind the algorithm. %
Assume $0 \leq y \leq x$ and define $r^*$ as $x/y$. One has
\begin{equation}
\label{defz*}
\sqrt{x^2+y^2} = x + \frac{y}{z^*},
\end{equation}
with
\begin{equation}
\label{computez*simplecase}
z^* =  r^* + \sqrt{1+(r^*)^2}.
\end{equation}
Two cases need be considered. If {$x \geq 2y$}, then the simple use of
(\ref{computez*simplecase}) suffices. An overflow may occur (typically
when computing $r^*$ or its square) but in such a case, $y$ is
negligible in front of $x$, so that $\sqrt{x^2+y^2} \approx x$ with very good accuracy, and one easily checks that it is the value returned by the algorithm. If $y \leq x \leq 2y$, more care is needed. The main idea is that when a variable $a$ is of the form $c + s$, where $c$ is a constant and $s$ is small, we retain more accuracy by representing it by the FP number nearest $s$ than by a FP approximation to $a$. Thus,  as $r^*=x/y$ is close to $1$, it can be represented with better accuracy as $1+r_2$, where $r_2$ is the FP number nearest
$
r_2^* := r^* - 1.
$
Furthermore, Lemma~\ref{lemma-sterbenz} implies that $\delta = x-y$ is computed exactly, so that $r_2$ is obtained through the FP division of $\delta$ by $y$.
Now, in order to express $z^*$ in terms of $r_2^*$, from (\ref{computez*simplecase}) a starting point is
\begin{equation}
\label{computez*debutsecondcase}
\left\{\begin{array}{lll}
z^* &:=& \sqrt{2+r_3^*} + r_2^*+1,\mbox{~with}\\
r_3^* &:= & (r_2^*)^2+2r_2^*.
\end{array}\right.
\end{equation}
Unfortunately, (\ref{computez*debutsecondcase}) cannot be used as is. When $r_3^*$ is small (i.e., when $y$ is close to $x$), much information on $r_3^*$ is lost in the FP addition $2+r_3^*$. A large part of this information can be retrieved using
\begin{equation}
\label{alternate-sqrt-2+r3}
\sqrt{2+r_3^*} = \frac{r_3^*}{\sqrt{2}+\sqrt{2+r_3^*}}+\sqrt{2}.
\end{equation}
More precisely, if the \emph{computed value} of  $2+r_3^*$ is  $2+r_3^*+\epsilon$, the influence of that error $\epsilon$ on the computed value of $\sqrt{2+r_3^*}$ is (at order $1$ in $\epsilon$)
\[
\sqrt{2+r_3^*+ \epsilon} - \sqrt{2+r_3^*} \approx \frac{\epsilon}{2 \sqrt{2+r_3^*}},
\]
whereas, using (\ref{alternate-sqrt-2+r3}), the influence of the error $\epsilon$ becomes
\[
\left(\frac{r_3^*}{\sqrt{2}+\sqrt{2+r_3^*+\epsilon}}+\sqrt{2}\right) - \sqrt{2+r_3^*} \approx - \frac{\epsilon r_3^*}{2 \sqrt{2+r_3^*}\left(\sqrt{2}+\sqrt{2+r_3^*}\right)^2},
\]
which is significantly smaller. So, instead of (\ref{computez*debutsecondcase}), the algorithm uses the formula
\begin{equation}
\label{z*final}
z^* = \frac{r_3^*}{\sqrt{2}+\sqrt{2+r_3^*}}+r_2^* + 1 + \sqrt{2}.
\end{equation}
To implement (\ref{z*final}) accurately, $1+\sqrt{2}$ is approximated by the unevaluated sum of two FP numbers $P_h$ and $P_\ell$, and the summation is performed small terms first.

The difficulty of the analysis comes from the number of steps of the
algorithm, which leads to a large number of variables in the polynomial
systems. This requires a very carefully exploitation of the special
structures of these systems. This is made easier by the fact that, due to  the
careful design of the algorithm, the partial derivatives of the
relative error with respect to the individual errors have signs
that are easily evaluated.

\subsection{Automatic Analyses}
This algorithm is too long for our current implementation. Still,
partial results are obtained automatically: a complete answer for
the first path ($\delta>y$) and $u\le 1/256$,
and an analysis of the linear term of the error for $\delta\le y$. In
order to obtain tighter bounds, this latter case is itself split into
the ranges $\alpha\in[1/2,2/3]$ and $\alpha\in[2/3,1]$. 
\paragraph{First path}
With input
\begin{Verbatim}
> Algo5firstpath:=[Input(x=0..2^16,alpha=1/2^16..1/2,_u=0..1/256),
> y=RN(alpha*x,0),r=RN(x/y),t=RN(1+r^2),s=RN(sqrt(t)),z=RN(r+s),
> z2=RN(y/z),rho=RN(x+z2)]:
> BoundRoundingError(Algo5firstpath);
\end{Verbatim}
Our program returns
\begin{multline*}
\color{blue}\left(\frac{157}{10}-\frac{32 \sqrt{5}}{5}\right) 
\textit{\_u} +\left(-\frac{347776}{5}-\frac{176400770811745024 
\sqrt{5}}{2155176452406911}\right.\\
\color{blue}\left.+\frac{27809906688 
\sqrt{88443606565122}\, \sqrt{5}}{8385900593023}\right) \textit{\_u}^{2}
\end{multline*}
whose numerical value is
\begin{Verbatim}
> evalf(%);
\end{Verbatim}
\[\color{blue}1.38916495 \textit{\_u} - 0.45335 \textit{\_u}^{2}\]
\paragraph{Second path}
Here are the linear parts for both subcases:
\begin{Verbatim}
> Algo5secondpath:=y=RN(alpha*x,0),delta=RN(x-y),r2=RN(delta/y),
> tr2=RN(2*r2,0),r3=RN(tr2+r2^2),r4=RN(2+r3),s2=RN(sqrt(r4)),
> sqrt2=RN(sqrt(2)),d=RN(sqrt2+s2,absolute(2*_u)),q=RN(r3/d),
> Ph=RN(1+sqrt2),Pl=RN(1+sqrt(2)-Ph,absolute(_u^2)),r5=RN(Pl+q),
> r6=RN(r5+r2),z=RN(Ph+r6),z2=RN(y/z),rho= RN(x+z2);
> split1:=[Input(x=1/2^16..2^16,alpha=1/2..2/3),Algo5secondpath]:
> split2:=[Input(x=1/2^16..2^16,alpha=2/3..1),Algo5secondpath]:
> BoundRoundingError(split1,steps="linear");
> BoundRoundingError(split2,steps="linear");
\end{Verbatim}
\begin{gather*}\color{blue}
\left(\frac{1019}{65}-\frac{46 \sqrt{13}}{13}-\frac{126 \sqrt{26}}
{65}+6 \sqrt{2}\right) \textit{\_u}\\
\color{blue}\left(-2+\frac{5 \sqrt{2}}{2}\right) \textit{\_u}
\end{gather*}
The precise value given for the error in \texttt{Pl} as a second
argument to
\texttt{RN} is explained in the proof of \cref{thm:error-Kahan} below.

Numerically, the values that have been computed are
approximately $1.5198\_u$ and $1.5355\_u$, so that the second one dominates.

\subsection{Result}
Our main result for this algorithm proves the linear term above and
gives a precise bound on the quadratic term.
\begin{theorem}\label{thm:error-Kahan}
For $p\ge 5$, the relative error of Algorithm~\ref{alg-Kahan} is
bounded by
\begin{align*}
\phi(u)&:=\frac{\left(1+\frac{1+u \left(1-2 u \right)}{1+\sqrt{2}-\frac{u \left(2+u\right) \left(1+2 u\right)^{2}}{\left(1+u \right)^{2}}}\right) \left(1+2 u\right)}{ \left(1+\frac{1}{1+\sqrt{2}}\right)\left(1+u \right)}-1,\\
&=\left(\frac{5 \sqrt{2}}{2}-2\right) u +\left(30-\frac{39 \sqrt{2}}{2}\right) u^{3}+\mathrm{O}\! \left(u^{4}\right),\\
&\le \left(\frac{5 \sqrt{2}}{2}-2\right) u+\frac{u^2}{12},\qquad 0\le u\le \frac1{32}.
\end{align*}
\end{theorem}

Here are examples of actually attained errors:

\begin{itemize}
   \item if $u = 2^{-24}$ (binary32 format), then error $1.4977 u$ is attained with $x = 12285049$ and $y = 11439491$;
   \item if $u = 2^{-53}$ (binary64 format), then error $1.4961 u$ is attained with $x = 6595357501251898$ and $y = 6135139757867044$.
\end{itemize}

We do not know if the bound given by  Theorem~\ref{thm:error-Kahan} is sharp. However, As $\frac{5 \sqrt{2}}{2}-2 \approx 1.5355$, the above examples show that there is not much room for improvement.
The proof of \cref{thm:error-Kahan} is given in \cref{sec:proof-kahan}.

\section{Discussion and Comparison}
\label{sec:discussion}

Table~\ref{tab-results} summarizes the error bounds  obtained for the
various algorithms considered in this article. To the best of our
knowledge, all our error bounds are new, and even the linear term for Algorithms~\ref{simplest-scaling}, \ref{Beebe-alg}, and~\ref{alg-Kahan} was unknown.

\begin{table}[tbh]
\begin{small}\begin{tabular}{|c|c|c|c|c|}
\hline
Algorithm & reference & error bound & status of bound \\
\hline\hline
\ref{NaiveHypot} & \begin{tabular}{c} straightforward \\
formula \end{tabular} & $2u-\frac{8}5(9-4\sqrt6)u^2\ (p\ge2)$& 
\begin{tabular}{c} asymptotically \\ optimal  \end{tabular}\\
\hline
\ref{simplest-scaling} &  \begin{tabular}{c} computer\\ 
arithmetic \\ folklore \end{tabular}  &$\frac{5}{2} u + \frac{3}{8}
u^2\ (p\ge2)$& \begin{tabular}{c} asymptotically \\ optimal  
\end{tabular}\\
\hline
\ref{Beebe-alg} & N.~Beebe~\cite{Beebe:2017}& $\frac{8}{5} u
+ \frac75 u^2\ (p\ge4)$&sharp \\
\hline
\ref{Borges-fused} & C.~Borges~\cite{Borges2020} & $u +
13.1u^2\ (p\ge5)$ &\begin{tabular}{c} asymptotically \\ optimal  
\end{tabular}\\
\hline
\ref{alg-Kahan} & W.~Kahan~\cite{Kahan1987} & $\left(\frac{5 \sqrt{2}}{2}-2\right) u+\frac{u^2}{12}\ (p \ge 5)$ & ? \\
\hline
\end{tabular}
\end{small}
\caption{Error bounds for Algorithms~\ref{NaiveHypot}, 
\ref{simplest-scaling}, \ref{Beebe-alg}, \ref{Borges-fused}, 
and~\ref{alg-Kahan}. More precise estimates for Algorithms~
\ref{NaiveHypot}, \ref{Borges-fused} are
given in \cref{th:naive-alg,thm-alg-Borges}.}
\label{tab-results}
\end{table}

\subsection{Comparison with Gappa}

We used version 1.3.5 of Melquiond's tool \texttt{gappa}%
\footnote{\url{https://gappa.gitlabpages.inria.fr}}.
As Gappa does not provide generic error bounds, we had to assume a given precision. Here, we give some results assuming binary32 arithmetic ($p = 24$).

\subsubsection{Algorithm~\ref{NaiveHypot}}

We fed Gappa with an input file that describes the algorithm\footnote{The input Gappa files are available on arXiv with this
article.}.
The returned result is
\begin{Verbatim}
  relerr in [0, 72412313008905457b-79 {1.19796e-07, 2^(-22.9929)}]
\end{Verbatim}
 which means that for binary32/single-precision arithmetic (i.e., $u = 2^{-24}$), Gappa finds a relative error bound
 $72412313008905457\times 2^{-79} \approx 2.00984u$ which is excellent, only very slightly above the Rump-Jeannerod bound $2u$. 
 
 Note that if we already know the error bound, we can ask Gappa to
 check it. In the input file, if we remove
 the absolute values in the description of \texttt{relerr}, replace
\texttt{     -> relerr in ?} by \texttt{-> relerr in  [-1b-23,1b-23]} and give the hint \texttt{ relerr \$ x in 64, y}, Gappa confirms that the bound $2^{-23}$ (i.e., $2u$) is correct.

\subsubsection{Algorithm~\ref{simplest-scaling}}

The result returned by Gappa for this algorithm is
\begin{Verbatim}
 relerr in [0, 792222648878942097b-82 {1.63828e-07, 2^(-22.5413)}]
\end{Verbatim}
 which means that for binary32/single-precision arithmetic (i.e., $u =
 2^{-24}$), Gappa finds a relative error bound $792222648878942097
 \times 2^{-82} \approx   2.749 u$, which is quite good but
 significantly larger than our bound $\frac{5}{2} u + \frac{3}{8} u^2$. In contrast to the case of Algorithm~\ref{NaiveHypot}, asking Gappa to confirm our bound instead of asking it to find one does not work.

\subsubsection{Algorithms~\ref{Beebe-alg} and~\ref{Borges-fused} }

With Algorithm~\ref{Beebe-alg} (which is essentially Algorithm~\ref{simplest-scaling} followed by a Newton-Raphson correction), Gappa returns an  error bound around twice the one it returns for Algorithm 2: it fails to ``see'' that lines 7 to 10 of the algorithm are a correction, and considers that these additional lines bring additional rounding errors. We could not find a ``hint'' that, provided to Gappa, would have improved the situation. The same phenomenon occurs with Algorithm~\ref{Borges-fused}: Gappa cannot see that the last lines of the algorithm are a Newton-Raphson correction.

\subsubsection{Algorithm~\ref{alg-Kahan}}

For the ``difficult'' path of the algorithm: $y < x < 2y$, Gappa
returns
\begin{Verbatim}
relerr in [0, 438344170044734879b-67 {0.00297034, 2^(-8.39516)}]
\end{Verbatim}
i.e., the computed error bound is $0.00297034 \approx 49834 u$: Algorithm~\ref{alg-Kahan} is clearly too complex to be  handled adequately.

\subsection{Comparison with Satire}

We used version 1.1 of the Satire tool%
\footnote{\url{https://github.com/arnabd88/Satire}}.
As Satire does not provide generic error bounds, we had
to assume a given precision. Here, we give some results assuming
binary32 arithmetic ($p = 24$)\footnote{The input Satire files are
available on arXiv with this
article.}.
Also, Satire computes absolute rather than relative error bounds.
Such bounds can also be obtained by our program, with the optional
argument \texttt{type="absolute"}. The variables~$x$ and~$y$ had to be
restricted away from~0 so as to avoid an infinite error being reported
due to a possible negative rounding assumed for  the argument of the square
root.

\begin{table}
\begin{center}
\renewcommand{\arraystretch}{1.2}
\sisetup{table-alignment-mode=format,table-number-alignment=center,table-format=1.3e1}
\begin{tabular}{cSSSS[table-format=1.1e1]}
\toprule
\multirow{2}*{Algorithm}
&\multicolumn{1}{c}{\multirow{2}*{Satire}}
&\multicolumn{1}{c}{bound deduced}
&\multicolumn{1}{c}{our code}
&\multicolumn{1}{c}{our code}\\[-1mm]
&
&\multicolumn{1}{c}{from \cref{tab-results}}
&\multicolumn{1}{c}{(linear part)}
&\multicolumn{1}{c}{(quadratic part)}\\
\midrule
1&1.658e-2&1.105e-2&1.105e-2&-4.9e-10\\
2&3.301e-2&1.382e-2&1.358e-2&6.6e-10\\
3&3.577e-2&8.839e-3&8.287e-3&2.5e-10\\
4&5.121e3 &5.525e-3&5.525e-3&2.4e-9\\
5&8.911e11&8.483e-2&3.058e-2&\multicolumn{1}{c}{n.a.}\\
\bottomrule
\end{tabular}
\renewcommand{\arraystretch}{1.0}
\end{center}
\caption{Absolute errors reported by Satire and by our code with
$p=24$ and 
$(x,y)\in[0,2^{16}]^2$.}
\label{table:satire}
\end{table}

A comparison of the results is given in \cref{table:satire}.
The higher quality of the bounds produced by our approach has to be
put in balance with the difficulty of obtaining them. Our program is
limited to small algorithms, while Satire can analyze programs with
hundreds of lines. The last ``n.a.'' corresponds to Algorithm~5, where
our code is currently unable to compute a quadratic bound on the
error. As a cross-check on the values reported here, one can use the
bounds on the relative error from \cref{tab-results} and multiply them
by the largest possible value. This gives an upper bound on the
largest absolute error, displayed in the 3rd column of the table.
The bounds computed directly on the absolute error are 
either identical or only slightly smaller, showing that in practice,
computing bounds on the relative error  is sufficient in many
cases.

\subsection{Conclusion}

Our approach makes it possible to obtain generic analytic error bounds for programs that implement functions that are considered to be
 ``basic building blocks'' of numerical computation. Because it is partially automatic, it limits the risk of human error and allows us to work with programs that are small but significantly larger than those that can be reasonably handled by paper-and-pencil computation. We obtain bounds that are often sharp, sometimes even asymptotically optimal, and, in the  tested cases, tighter than those provided by the other existing tools. However, Satire can handle much larger programs, and Gappa has the valuable ability to provide formal proofs of the bounds it calculates.

An important problem illustrated by this work is the large place taken
by computations in the proofs. Thus in this area, more than in others,
we feel that an increasingly important role will have to be taken
by computer-aided proofs. This is part of a more general trend that
affects a large part of mathematics~\cite{Granville2024}. Currently,
certified computer algebra is still a long-term goal, despite progress
being made in that direction\footnote{See for instance 
\url{https://fresco.gitlabpages.inria.fr}.}. Thus we have provided both
pencil-proofs and link to Maple worksheets, but this will not be an
option for longer proofs.


\pagebreak
\appendix
\section{Proof of Theorem \ref{th-sharp-alg-beebe}}
\label{app-proof-beebe}
\subsection{Linear term and interval splitting}
The Taylor expansion of~$R$ in \cref{lemma:rel-err-3} shows the importance of a good upper bound on
the error on~$r$. If one uses the
underlined
equation in \cref{eq:algo3-part1},
then $r^2=\alpha^2(1+u\epsilon_r)^2$;
the linear
term is
maximized when $\epsilon_\rho,\epsilon_t,\epsilon_r$ are all equal to~1, where
it becomes
\[\frac{2\alpha^2+1}{2(1+\alpha^2)}+1=2-\frac{1}{2(1+\alpha^2)}.
\]
This bound on the linear term is increasing with respect to~$\alpha$.
If $\alpha$
ranges from 0 to~1, then the maximum is reached at
$\alpha=1$, where the value
is~$7/4$. This is the
result of our first automatic analysis.
If $\alpha$ ranges from 0 to $1/2$ only, then the bound reached
at~$\alpha=1/2$ is $8/5$, the result of our second automatic analysis.

Finally, if $\alpha\ge 1/2$, then so is its rounded value~$r$ and one
can use the \emph{absolute} error bound
\begin{equation}\label{eq:algo3-abs-err-r}
r=\alpha+\frac{u\epsilon_r}2,\quad |\epsilon_r|\le 1.
\end{equation}
The linear term above becomes
\[\frac{\alpha+1}{2(1+\alpha^2)}+1,\]
which is decreasing with~$\alpha$ for $\alpha\ge1/2$. Thus the
maximum is reached at $\alpha=1/2$, where again $8/5$ is obtained.

Actually, using an absolute error bound is beneficial for $\alpha\le 1/2$ as
well.  This type of variation on
the analysis can also be requested of our implementation, by adding the keyword
\texttt{absolute} as second argument to \texttt{RN}:
\begin{Verbatim}
> Algo3_1:=[Input(x=0..2^16,alpha=0..1/2),y=RN(alpha*x,0),
> r=RN(y/x,absolute),t=RN(1+r^2),s=RN(sqrt(t)),
> epsilon=RN(t-s^2,0),c = RN(epsilon/2/s,absolute(_u^2/2)),
> nu = RN(x*c),rho=RN(nu+x*s)]: 
> BoundRoundingError(Algo3_1,steps="linear");
\end{Verbatim}
\[
{\color{blue}\left(\frac{5}{4}+\frac{\sqrt{5}}{8}\right) \textit{\_u}}
\]
(The optional argument \texttt{steps="linear"} specifies that only
the linear term of the error bound is computed.)

Indeed, for $\alpha\le 1/2$, the absolute error bound becomes
\[r=\alpha+\frac{u\epsilon_r}4,\quad |\epsilon_r|\le 1.\]
While this bound is not as good as the relative error bound for $\alpha<1/4$,
it allows one to bound the linear term of \cref{eq:defR-algo3} by
\[\frac{\alpha/2+1}{2(1+\alpha^2)}+1,\]
whose maximum~---~reached at $\alpha=\sqrt{5}-2\simeq0.236$~---~equals $(10+
\sqrt{5})/8\simeq 1.53$, which is smaller than $8/5$. Thus the maximal value of
the linear part is reached as $\alpha$ approaches $1/2$ from the right.

\subsection{Bounds on partial derivatives}
In view of the previous discussion, we replace the underlined bound in 
\cref{eq:algo3-part1} by
\begin{equation}\label{eq:algo3-r}
r=\alpha+u\epsilon_r,\quad|\epsilon_r|\le\begin{cases}
\frac14,\quad&\text{if $\alpha\le1/2$,}\\ 
\frac12,\quad&\text{if $\alpha>1/2$.}
\end{cases}\end{equation}

The analysis of the linear term shows that in a neighborhood of $u=0$, the
relative error $R$, seen as a function of $u$,
$\alpha$ and the $\epsilon_i$, is increasing with
$\epsilon_\rho,\epsilon_t,\epsilon_r$. We
first show that these properties hold for all $u\in[0,1/2]$ and that
the derivative wrt $\epsilon_c$ also has constant sign.
\begin{lemma}\label{lemma:algo3-derivatives} For $u\in[0,1/2]$, 
\[\frac{\partial R}{\partial\epsilon_\rho}\ge0,\quad 
\frac{\partial R}{\partial\epsilon_t}\ge0,\quad 
\frac{\partial R}{\partial\epsilon_r}\ge0,\quad\frac{\partial R}
{\partial\epsilon_c}\ge0.
\]
\end{lemma}
\begin{proof}
Simple inequalities that can be used in the analysis are
\begin{equation}\label{eq:algo3-simple-ineqs}
|\epsilon_i|\le 1\quad i\in\{r,s,c,\nu,\rho\},\quad 
r\in[0,1],\quad t\in[1,2],\quad s\in[1,2].
\end{equation}
For $u\in[0,1/2]$, we first consider the factors of $R+1$ in 
\cref{eq:defR-algo3}:
\begin{gather*}
\left|\frac{u^2}{2\sqrt{t}}
((\epsilon_c+\epsilon_s^2/s)(1+u\epsilon_\nu)-2\epsilon_s\epsilon_\nu)\right|
\le\frac{u^2}2(2(1+u)+2)=u^2(2+u)<1;\\
\left|\frac{r^2-\alpha^2}{1+\alpha^2}\right|=
\left|\frac{u\epsilon_r
(2\alpha+u\epsilon_r)}{1+\alpha^2}\right|\le \frac{u(2+u/2)}2<1.
\end{gather*}
This implies that all the factors of $R+1$ in \cref{eq:defR-algo3}
are positive. Also, it follows that the
derivative of~$R$ with respect to~$\epsilon_\rho$ is positive.

Using \crefrange{eq:algo3-part1}{eq:algo3-part2}, it can be checked that
\[\frac{\partial R}{\partial\epsilon_t}=\frac{(R+1)ux(2
s^2-u^2\epsilon_s^2(1+u\epsilon_\nu))}{4s^2\sqrt{t}
(\nu+xs)}\]
and the nonnegativity of the right-hand side for $u\in[0,1/2]$ is
clear since $s\ge1$. The same conclusion follows for $\epsilon_r$ in
view of
\[\frac{\partial R}{\partial\epsilon_r}=2r\frac{\partial R}
{\partial\epsilon_t},\]
which, again, is checked with \crefrange{eq:algo3-part1}{eq:algo3-part2}.
Finally, the case of the derivative with respect to $\epsilon_c$ follows from
\[\frac{\partial R}{\partial\epsilon_c}=(R+1)\frac{u^2x(1+u\epsilon_\nu)}{2
(xs+\nu)}.\qedhere\]
\end{proof}

\subsection{Bounds on the relative error}

\subsubsection{First part: upper bound for
\texorpdfstring{$\alpha<1/2$}{alpha<1/2} and lower bound for
 all \texorpdfstring{$\alpha$}{alpha}}
For these cases, we prove bounds on the absolute value of~$R$ that are smaller than
$8u/5+u^2$ for $u\in[0,1/16]$.

From \cref{eq:algo3-simple-rho}, it follows that
\[R=\frac{1}{\sqrt{1+\alpha^2}}\left(\sqrt{t}+\frac{u^2}2
((\epsilon_c+\epsilon_s^2/s)(1+u\epsilon_\nu)-2\epsilon_s\epsilon_\nu)
\right)
(1+u\epsilon_\rho)-1.\]
Upper (resp. lower) bounds on the factor of~$u^2/2$ are reached at $
(\epsilon_c,\epsilon_s,\epsilon_\nu)=(1,-1,1)$ (resp. $(-1,1,1))$).
Thus,
\begin{multline*}
\frac{1}{\sqrt{1+\alpha^2}}\left(\sqrt{t}-{u^2}\frac{3s-1}{2s}-u^3
\frac{s-1}{2s}\right)
(1+u\epsilon_\rho)\\ \le R+1
\le\frac{1}{\sqrt{1+\alpha^2}}\left(\sqrt{t}+
{u^2}\frac{3s+1}{2s}+u^3\frac{s+1}{2s}
\right)
(1+u\epsilon_\rho)
\end{multline*}
The lower bound is lower bounded by its value at $s=2$, while the upper
bound is upper bounded by
its value at~$s=1$, whence
\begin{multline*}
\frac{1}{\sqrt{1+\alpha^2}}\left(\sqrt{t}-\frac{5u^2}{4}-
\frac{u^3}{4}\right)
\left(1-\frac{u}{1+u}\right)\\ \le R+1
\le\frac{1}{\sqrt{1+\alpha^2}}\left(\sqrt{t}+
2{u^2}+u^3
\right)
\left(1+\frac{u}{1+u}\right)
\end{multline*}
Next, $t=1+u\epsilon_t+(\alpha+u\epsilon_r)^2$ is minimal for $\epsilon_t=-1$
and $\epsilon_r=-\epsilon_r^{\max}$ and maximal for the opposite of these
values, giving
\begin{multline*}
\frac{1}{\sqrt{1+\alpha^2}}\left(\sqrt{1-u+(\alpha-u\epsilon_r^{\max})^2}-
\frac{5u^2}{4}-
\frac{u^3}{4}\right)
\left(1-\frac{u}{1+u}\right)\\ \le R+1
\le\frac{1}{\sqrt{1+\alpha^2}}\left(\sqrt{1+u+(\alpha+u\epsilon_r^{\max})^2}+
2{u^2}+u^3
\right)
\left(1+\frac{u}{1+u}\right)
\end{multline*}
We show the steps of the proof that 
\[R\ge -\frac85u\quad\text{for}\quad\alpha\in[0,1/2],\quad\epsilon_r^{
\max{}}=1/4,\quad
u\in
[0,1/4].\] The proofs of the following two other inequalities,
\begin{align*}
R\ge -\frac85u&\quad\text{for}\quad\alpha\in[1/2,1],\quad\epsilon_r^{
\max{}}=1/2,\quad u\in
[0,1/4],\\
R\le \frac85u+u^2&\quad\text{for}\quad\alpha\in
[0,1/2],\quad\epsilon_r^{
\max{}}=1/4,\quad u\in[0,1/16],
\end{align*}
follow the same template.

The function 
\[\frac{1}{\sqrt{1+\alpha^2}}\left(\sqrt{1-u+(\alpha-u/4)^2}-
\frac{5u^2}{4}-
\frac{u^3}{4}\right)
\left(1-\frac{u}{1+u}\right)-1+\frac85u\]
is positive for $u\rightarrow0$, $\alpha\in[0,1/2]$. Eliminating both square
roots shows that it does not vanish unless the following polynomial does:
\begin{multline*}
E(u,\alpha)=1600 \left(12 \alpha^{2}-5 \alpha +2\right)^{2}+80 \left(1136
\alpha^{2}+1161\right) \left(12 \alpha^{2}-5 \alpha +2\right) u\\ 
 +\left(553216 \alpha^{4}+307200 \alpha^{3}+777632 \alpha^{2}+307200 \alpha
 +225041\right) u^{2}\\ 
 +\left(-2727936 \alpha^{4}+409600 \alpha^{3}-4475072 \alpha^{2}+659600 \alpha
 -1247136\right) u^{3}\\ 
 +\left(-1736704 \alpha^{4}-1904608 \alpha^{2}+100000 \alpha +846\right) u^
 {4}\\ 
 +\left(1572864 \alpha^{4}+2777728 \alpha^{2}+10000 \alpha +1212364\right) u^
 {5}\\ 
 +\left(1048576 \alpha^{4}+489952 \alpha^{2}-169249\right) u^{6}\\ 
 +\left
 (-550400 \alpha^{2}-237900\right) u^{7}\\ 
 +\left(-51200 \alpha^{2}+42550\right)
 u^{8}+12500 u^{9}+625 u^{10}.
\end{multline*}
The conclusion follows from the fact that $E(u,\alpha)>0$ for $\alpha\in
[0,1/2]$, $u\in[0,1/4]$. This can be checked by minimizing each of the
coefficients over the interval $\alpha\in[0,1/2]$ and then evaluating at
$u\in[0,1/4]$ by interval arithmetic.

\subsubsection{Last part: upper bound for
\texorpdfstring{$1/2\le\alpha\le1$}{1/2≤alpha≤1}}
In order to obtain a tight bound on the relative error in that case, we first
consider the derivative with respect to $\epsilon_\nu$ and $\epsilon_s$.
\paragraph{Derivative with respect to $\epsilon_\nu$.}
Similar manipulations as above give
\[\frac{\partial R}{\partial\epsilon_\nu}=(R+1)\frac{ucx}{xs+\nu},
\]
whose sign is given by that of~$c$. Thus the extremal value are reached at
$\epsilon_\nu=\pm1/(1+u)$, except when $c=0$.

The case $c=0$ can be treated separately. When $c=0$, then $\epsilon=\nu=0$, $s=
\sqrt{t}$ and $R$ is bounded by~$8u/5$ for $\alpha\in[1/2,1],u\in
[0,1/2],|\epsilon_r|\le1/2$:
\begin{align*}R=\frac{\sqrt{t}}{\sqrt{1+\alpha^2}}-1
&=\sqrt{1+u\frac{2\alpha\epsilon_r+\epsilon_t+u\epsilon_r^2}{1+\alpha^2}}-1\\ 
&\le \sqrt{1+u\frac{1+\alpha+u/4}{1+\alpha^2}}-1
\le\sqrt{1+u(1+u/5)}-1
\le\frac85u.
\end{align*}
\paragraph{Derivative with respect to $\epsilon_s$.}
It can be checked with \crefrange{eq:algo3-part1}{eq:algo3-part2} that
\[\frac{\partial R}{\partial\epsilon_s}=\frac{u^2(1+u\epsilon_\rho)}
{2s^2\sqrt{1+\alpha^2}}\left(-2t\epsilon_\nu+2\sqrt{t}\epsilon_s+u\epsilon_s
(\epsilon_s-2\sqrt{t}\epsilon_\nu)-u^2\epsilon_\nu\epsilon_s^2\right).\]
When $\epsilon_\nu=1/(1+u)$ and $\alpha\in[1/2,1]$, the last factor
reaches its maximum for $(u,\epsilon_s,t)\in[0,1/4]\times
[-1,1]\times[5/4,2]$ at $t=5/4$, $u=1/4$, $\epsilon_s=1$, where it equals
$4/\sqrt5-9/5<0$,
making the derivative negative. Similarly, when
$\epsilon_\nu=-1/(1+u)$ and $\alpha\in[1/2,1]$, the last factor
reaches its minimum for $(u,\epsilon_s,t)\in[0,1/16]\times[-1,1]\times[5/4,2]$
at $t=5/4$, $u=1/16$, $\epsilon_s=-1$, where it equals
$(329-144\sqrt5)/136>0$,
making the derivative positive.

Thus, for $\alpha\in[1/2,1]$ and $u\in[0,1/16]$, at the extremal values
$\epsilon_\nu=\pm1/(1+u)$, the value of~$R$ is
maximized when $\epsilon_s=\mp1$.

\paragraph{Derivative with respect to $\alpha$}
The derivative with respect to~$\alpha$ factors as
\begin{multline*}
\frac{\partial R}{\partial\alpha}=
\frac{u(1+u\epsilon_\rho)}
{2s^2\sqrt{t}(1+\alpha^2)^{3/2}}
\left(
 -2 \left(\alpha^{2} \epsilon_{r}+\alpha  \epsilon_{t}-\epsilon_{r}\right) s^
 {2}\right.\\
+\left(2 \sqrt{t}\, \alpha  \,s^{2} \epsilon_{\nu} \epsilon_{s}-\sqrt{t}\,
\alpha  \,s^{2} \epsilon_{c}-\sqrt{t}\, \alpha  s \epsilon_{s}^{2}-\alpha^{2} r \epsilon_{s}^{2}-2 \alpha  \epsilon_{r}^{2} s^{2}-r \epsilon_{s}^{2}\right) u\\
\left.-\epsilon_{\nu} \left(\sqrt{t}\, \alpha  \,s^{2} \epsilon_{c}+\sqrt{t}\,
\alpha  s \epsilon_{s}^{2}+\alpha^{2} r \epsilon_{s}^{2}+r \epsilon_{s}^{2}\right) u^{2}
\right).
\end{multline*}

By \cref{lemma:algo3-derivatives} and the evaluation of the derivatives above,
for $\alpha\in[1/2,1]$ and $u\in
[0,1/16]$, $|R|$ is upper
bounded by the maximum of its values at the points
\[
(\epsilon_\rho,\epsilon_t,\epsilon_r,\epsilon_c,\epsilon_\nu,\epsilon_s)=
\left(\frac1{1+u},1,\frac12,1,\frac{\pm1}{1+u},\mp1\right):=\pi_\pm.\]
At these values, the last factor in $\partial R/\partial\alpha$ is upper bounded
by
\[-2\left(\frac18\right)-\left(\frac32+\frac12+\frac18+\frac14+\frac12\right)u
+(4\sqrt2+2\sqrt2+1+1)u^2,\]
which is negative for $u\in[0,1/16]$ making the derivative with respect
to~$\alpha$ negative, so that the maximum is reached at $\alpha=1/2$.

\paragraph{Final bound}
By \cref{lemma:algo3-derivatives} and the evaluation of the derivatives above,
the maximal value of~$R$ is reached at $\alpha=1/2$ and $
(\epsilon_\rho,\epsilon_t,\epsilon_r,\epsilon_c,\epsilon_\nu,\epsilon_s)=\pi_\pm$, where
\[R=\frac{1}{\sqrt{1+1/4}}\left(\sqrt{t}+\frac{u^2}2
\left((1+1/s)\left(1\pm\frac u{1+u}\right)+\frac2{1+u}\right)
\right)
\left(1+\frac u{1+u}\right)-1,\]
with
\[s=\sqrt{t}\mp u,\quad t=1+u+\left(\frac12+\frac u2\right)^2.\]
It follows that the maximal value is reached at~$\pi_+$.

The approach for the last step is as in the analysis of Algorithm~2.
At $\pi_+$,
the quantity~$y=(R-\frac85u)/u^2$ to be maximized has Taylor expansion
\[y=\frac{3}{\sqrt 5}-\frac2{25}+O(u)\]
and satisfies the quartic
equation $E(y,u)=0$ with
\begin{multline*}
E(y,u)=
625 u^{4} \left(3 u^{2}-6 u -5\right)^{2} \left(1+u \right)^{6} y^{4}\\ 
+500 u^
{2} \left(8 u +5\right) \left(3 u^{2}-6 u -5\right)^{2} \left(1+u \right)^{6} y^
{3}\\ 
-50 \left(720 u^{11}-2640 u^{10}-6788 u^{9}+9528 u^{8}+44735 u^{7}+54117
u^{6}+5941 u^{5}\right.\\ 
\left.-55797 u^{4}-67930 u^{3}-37950 u^{2}-10750 u -1250\right)
\left(1+u \right)^{3} y^{2}\\ 
-20 \left(8 u +5\right) \left(720 u^{9}-2640 u^
{8}-5636 u^{7}+9816 u^{6}+34145 u^{5}\right. \\ 
\left.+38979 u^{4}+22767 u^{3}+6481 u^
{2}+460 u -100\right) \left(1+u \right)^{3} y\\ 
 +57600 u^{14}-192000 u^
{13}-940160 u^{12}+781440 u^{11}+7808624 u^{10}\\ 
+11779648 u^{9}-5174360 u^
{8}-47142584 u^{7}-88187555 u^{6}-95250006 u^{5}\\ 
-66769667 u^{4}-30792720 u^
{3}-9017844 u^{2}-1519640 u -112100
\end{multline*}
The function~$y$ is increasing in the range $u\in[0,1]$, as can be seen by
checking that the resultant of this polynomial with its derivative with respect
to~$u$ does not vanish there. Thus the bound $Q$ on the quadratic term for $u\in
[0,u_{\max{}}]$ with
$u_{\max{}}\le1/16$ is maximal
at~$u=u_{\max{}}$, where it has the values given in \cref{table:quad-algo3}.
(Explicit, but not very useful, expressions in terms of radicals are available
for these bounds.)

\section{Proof of Theorem~\ref{thm-alg-Borges}}
\label{app-proof-thm-Borges}
\subsection{Step-by-step analysis}
From the properties of Algorithms Fast2Sum and Fast2Mult and \cref{lemma-knuth-dekker,lemma-boldodaumasSQRT}, we obtain
the following relations (as before, \(u\le 1/4\)):
\begin{align}
y=\alpha x,&\notag\\
s_x^h=x^2(1+u\epsilon_{s_x^h}),&\quad |\epsilon_{s_x^h}|\le\frac1
{1+u},\label{eq:sxh}\\
s_x^h+s_x^\ell=x^2,&\label{eq:sxl}\\
s_y^h=y^2(1+u\epsilon_{s_y^h}),&\quad |\epsilon_{s_y^h}|\le\frac1
{1+u},\label{eq:syh}\\
s_y^h+s_y^\ell=y^2,&\label{eq:syl}\\
\sigma_h=(s_x^h+s_y^h)(1+u\epsilon_{\sigma_h}),
&\quad|\epsilon_{\sigma_h}|\le\frac1{1+u},\label{eq:sigmah}\\
\sigma_h+\sigma_\ell=s_x^h+s_y^h,&\label{eq:sigmal}\\
s=\sqrt{\sigma_h}(1+u\epsilon_s),&\quad|\epsilon_s|\le \frac1u
(1-\frac1{\sqrt{1+2u}})=\frac{2}{\sqrt{1+2u}(1+\sqrt{1+2u})},
\label{eq:s}\\
s^2+\delta_s=\sigma_h,&\label{eq:deltas}\\
\tau_1=(s_x^\ell+s_y^\ell)(1+u\epsilon_{\tau_1}),
&\quad|\epsilon_{\tau_1}|\le \frac1{1+u},\label{eq:tau1}\\
\tau_2=(\delta_s+\sigma_\ell)(1+u\epsilon_{\tau_2}),
&\quad|\epsilon_{\tau_2}|\le \frac1{1+u},\label{eq:tau2}\\
\tau=(\tau_1+\tau_2)(1+u\epsilon_\tau),&\quad|\epsilon_\tau|\le \frac1
{1+u},\label{eq:tau}\\
c=\frac{\tau}s(1+u\epsilon_c),&\quad|\epsilon_c|\le 1-2u,
\label{eq:c}\\
\rho=\left(s+\frac12c\right)
(1+u\epsilon_\rho),&\quad|\epsilon_\rho|\le \frac1{1+u}.\label{eq:rho}
\end{align}

For the analysis, it is important to isolate terms of order \(u\), thus
we introduce variables \(\theta_v\) so that
\[v=u\theta_v,\qquad v\in\{s_x^\ell,s_y^\ell,\sigma_\ell,\delta_s,\tau_1,\tau_2,\tau,c\}.\]
From the equations above, it follows that these are given explicitly by
\begin{align}
\theta_{s_x^\ell}&=-x^2\epsilon_{s_x^h},&\quad
\theta_{s_y^\ell}&=-y^2\epsilon_{s_x^h},\label{eq:theta1}\\
\theta_{\sigma_\ell}&=-(s_x^h+s_y^h)\epsilon_{\sigma_h},&\quad
\theta_{\delta_s}&=-\sigma_h\epsilon_s(2-u\epsilon_s),\label{eq:theta2}\\
\theta_{\tau_1}&=(\theta_{s_x^\ell}+\theta_{s_y^\ell})(1+u\epsilon_
{\tau_1}),&\quad
\theta_{\tau_2}&=(\theta_{\delta_s}+\theta_{\sigma_\ell})(1+u\epsilon_
{\tau_2}),\label{eq:theta3}\\
\theta_\tau&=(\theta_{\tau_1}+\theta_{\tau_2})(1+u\epsilon_\tau),&\quad
\theta_c&=\frac{\theta_\tau}s(1+u\epsilon_c).\label{eq:theta4}
\end{align}

\subsection{Relative error and its linear
term}\label{relative-error--its-linear-term}

The computation proceeds by using the \emph{equalities} defining the
successive
variables in such a way that the leading coefficient in \(u\) is made
apparent, and the error term is given an explicit expression in terms of
the variables \(\epsilon_i\). This will then be followed by a
computation of bounds using the known \emph{inequalities} on these
quantities.

The starting point is
\begin{align}
x^2+y^2&=s_x^h+s_y^h+s_x^\ell+s_y^\ell,&&\text{from \cref{eq:sxl,eq:syl}}\\
&=\sigma_h+\sigma_\ell+s_x^\ell+s_y^\ell,&&\text{from \cref{eq:sigmal}}\\
&=s^2+\delta_s+\sigma_\ell+s_x^\ell+s_y^\ell,&&\text{from 
\cref{eq:deltas}}\\
&=s^2\left(1+\frac{\delta_s+\sigma_\ell+s_x^\ell+s_y^\ell}
{s^2}\right).&&
\end{align}
Then, taking square roots,
\begin{align}
\sqrt{x^2+y^2}&=s\sqrt{1+\frac{\delta_s+\sigma_\ell+s_x^\ell+s_y^\ell}
{s^2}},\notag\\
&=s\sqrt{1+\frac1{s^2}(\frac{\tau_1}{1+u\epsilon_{\tau_1}}+
\frac{\tau_2}{1+u\epsilon_{\tau_2}})},\qquad\qquad\text{from 
\cref{eq:tau1,eq:tau2}}\notag\\
&=s\sqrt{1+\frac{\tau_1+\tau_2}{s^2}-\frac u{s^2}(\frac{\epsilon_
{\tau_1}\tau_1}{1+u\epsilon_{\tau_1}}+\frac{\epsilon_{\tau_2}\tau_2}
{1+u\epsilon_{\tau_2}})},\notag\\
&=s\sqrt{1+\frac\tau{s^2}-\frac u{s^2}(\frac{\epsilon_\tau\tau}
{1+u\epsilon_\tau}+\frac{\epsilon_{\tau_1}\tau_1}{1+u\epsilon_
{\tau_1}}+\frac{\epsilon_{\tau_2}\tau_2}{1+u\epsilon_{\tau_2}})},\quad
\text{from \cref{eq:tau}}\notag\\
&=s\sqrt{1+\frac{\tau}{s^2}-u^2A},\qquad\qquad\qquad\qquad\qquad
\text{from
\cref{eq:theta3,eq:theta4}}\notag\\
&\text{with}\quad A=\frac{1}{s^2}(
\frac{\epsilon_\tau\theta_\tau}{1+u\epsilon_\tau}+\frac{\epsilon_
{\tau_1}\theta_{\tau_1}}{1+u\epsilon_{\tau_1}}+\frac{\epsilon_
{\tau_2}\theta_{\tau_2}}{1+u\epsilon_{\tau_2}}),\label{eq:A}\\
&=s\left(1+\frac{\tau}{2s^2}-u^2B\right),\notag\\
&\text{with}\quad B=\left(1+u
\frac{\theta_\tau}{2s^2}-\sqrt{1+u\frac{\theta_\tau}
{s^2}-u^2A}\right)/u^2,\label{eq:B}\\
&=s+\frac c2-u^2\frac{\epsilon_c\theta_c}{2
(1+u\epsilon_c)}-u^2sB,\qquad\qquad\text{from
\cref{eq:c,eq:theta4}}\notag\\
&=\frac{\rho}{1+u\epsilon_\rho}-u^2\frac{\epsilon_c\theta_\tau}
{2s}-u^2sB,\qquad\qquad\quad\text{\ from
\cref{eq:rho}}\notag\\
&=\rho\left(\frac{1}{1+u\epsilon_\rho}-u^2C\right),\qquad C=
\frac{1}\rho\left(\frac{\epsilon_c\theta_\tau}{2s}+sB\right)
\label{eq:def-C}.
\end{align}

In summary, we have obtained the following.
\begin{lemma}\label{lemma-Borges} The relative error of 
Algorithm~\ref{Borges-fused} is
\begin{align*}
R:=\frac\rho{\sqrt{x^2+y^2}}-1=\frac1{\frac1{1+u\epsilon_\rho}-{u^2}C}-1=u\epsilon_\rho+O(u^2),
\end{align*}
with $C$ defined in \cref{eq:def-C}, in terms of variables $A,B$ and
\(\theta_v\) defined in 
\cref{eq:A,eq:B,eq:theta1,eq:theta2,eq:theta3,eq:theta4}. Moreover, in
this
formula,
\[|\epsilon_s|\le 1-\frac1{\sqrt{1+2u}},\qquad |\epsilon_c|\le 1-2u,\]
and
\(|\epsilon_{s_x^h}|,|\epsilon_{s_y^h}|,|\epsilon_{\sigma_h}|,|\epsilon_{\tau_1}|,|\epsilon_{\tau_2}|,|\epsilon_\tau|,|\epsilon_\rho|\)
are bounded by \(1/(1+u)\).
\end{lemma}

\subsection{Bounds}\label{bounds}
Our implementation claims that the maximal value of the bound
from \cref{lemma-Borges} is reached with the extreme values
\begin{gather*}
\epsilon_{s_x^h}=\epsilon_{s_y^h}=\epsilon_{\sigma_h}=-\frac1
{1+u},\quad
\epsilon_{\tau_1}=\epsilon_
{\tau_2}=\epsilon_\tau=\epsilon_\rho=\frac1{1+u},\\
\epsilon_c=1-2u,\quad
\epsilon_s=\frac1u\left(\frac1{\sqrt{1+2u}}-1\right).
\end{gather*}
We do not
have a human-readable proof of this fact, which has \cref{conj-beebe}
as a direct consequence.

\subsection{Proof of Theorem \ref{thm-alg-Borges}}
The proof consists in propagating bounds from variable to variable,
trying not to lose too much correlation between variables along
the way.
\begin{small}
\begin{align}
\frac1{1+u}x^2&\le s_x^h\le\frac{1+2u}{1+u}x^2,
\qquad\text{from \cref{eq:sxh}}\label{ineq:sxh}\\
\frac1{1+u}y^2&\le s_y^h\le \frac{1+2u}{1+u}y^2,
\qquad\text{from \cref{eq:syh}}\label{ineq:syh}\\
\frac1{(1+u)^2}&\le\frac{\sigma_h}{x^2+y^2}\le \frac{(1+2u)^2}{(1+u)^2},
\ \text{from \cref{eq:sigmah,ineq:sxh,ineq:syh}}\label{ineq:sigmah}\\
|\theta_{s_x^\ell}|&\le x^2\frac{1}{1+u},
\qquad\text{from \cref{eq:sxh,eq:theta1}}\label{ineq:thetasxl}\\
|\theta_{s_y^\ell}|&\le y^2\frac{1}{1+u},
\qquad\text{from \cref{eq:syh,eq:theta1}}\label{ineq:thetasyl}\\
|\theta_{\sigma_\ell}|&\le(x^2+y^2)\frac{1+2u}{(1+u)^2},
\quad\text{from \cref{eq:theta2,ineq:sxh,ineq:syh}}
\label{ineq:thetasigmal}\\
|\epsilon_s(2-u\epsilon_s)|&\le2\frac{2+3u-2\sqrt{1+2u}}{u(1+2u)},
\qquad\text{from \cref{eq:s}}\label{ineq:epsilons}\\
|\theta_{\delta_s}|&\le (x^2+y^2)\frac{(1+2u)}{(1+u)^2}\frac{4+6u-4\sqrt{1+2u}}{u},
\ \text{from \cref{eq:theta2,ineq:epsilons}}
\label{ineq:thetadeltas}\\
|\theta_{\tau_1}|&\le(x^2+y^2)\frac{1+2u}{(1+u)^2},
\quad\text{from 
\cref{eq:tau1,eq:theta3,ineq:thetasxl,ineq:thetasyl}}
\label{ineq:thetatau1}\\
|\theta_{\tau_2}|&\le(x^2+y^2)\frac{(1+2u)^2}{(1+u)^3}\frac{4+7u-\sqrt{1+2u}}u,
\ \text{from 
\cref{eq:tau2,eq:theta3,ineq:thetasigmal,ineq:thetadeltas}}
\label{ineq:thetatau2}\\
|\theta_\tau|&\le(x^2+y^2)\frac{(1+2u)^2}{(1+u)^4}\phi(u),
\ \text{from 
\cref{eq:tau,eq:theta4,ineq:thetatau1,ineq:thetatau2}}
\label{ineq:thetatau}\\
&\qquad\text{with }\phi(u):=\frac{(2+3u)(2+5u)-4(1+2u)^{3/2}}{u},\\
s&\ge\frac{\sqrt{x^2+y^2}}{(1+u)\sqrt{1+2u}}
\qquad\text{from \cref{eq:s,ineq:sigmah}}
\label{ineq:s-low}\\
s&\le\frac{\sqrt{x^2+y^2}(1+2u)\left(2-\frac1{\sqrt{1+2u}}\right)}{
(1+u)}
\qquad\text{from \cref{eq:s,ineq:sigmah}}
\label{ineq:s-up}\\
|A|&\le\frac1{s^2}(|\theta_\tau|+|\theta_{\tau_1}|+|\theta_
{\tau_2}|),\ \text{from \cref{eq:A},\eqref{ineq:sigmah} and }
\frac{\epsilon_{\tau_i}}{1-u\epsilon_{\tau_i}}=1,\label{ineq:A}\\
&\le\frac{(1+2u)^2(2+3u)}{(1+u)^2}\phi(u),
\qquad\text{from \cref{ineq:thetatau1,ineq:thetatau2,ineq:thetatau}}\\
c
&\ge-\sqrt{x^2+y^2}u\frac{(1-u)(1+2u)^{7/2}}{(1+u)^3}\phi(u),
\ \text{from \cref{eq:c,eq:theta4,ineq:thetatau,ineq:s-low}}
\label{ineq:c}\\
\rho&\ge\frac{\sqrt{x^2+y^2}}{1+u}\left(\frac1{(1+u)\sqrt{1+2u}}-\frac
u2\frac{(1-u)(1+2u)^{7/2}}{(1+u)^3}\phi(u)\right),\notag\\
&\qquad\qquad \text{from \cref{eq:rho,ineq:s-up,ineq:s-low,ineq:c}}
\label{ineq:rho}\\
&\ge\frac{\sqrt{x^2+y^2}}{(1+u)^2\sqrt{1+2u}}\left(1-\frac u2\frac{(1-u)(1+2u)^{4}}{(1+u)^2}\phi(u)\right),\\
|B|&\le\left(1-\frac u2\frac{(1+2u)^3}{(1+u)^2}\phi(u)-\sqrt{1-u\frac{(1+2u)^3}{(1+u)^2}\phi(u)-u^2\frac{(1+2u)^2(2+3u)}{(1+u)^2}\phi(u)}\right)/u^2,\label{ineq:B}\\
&\qquad\text{from \cref{eq:B,ineq:thetatau,ineq:s-low,ineq:A} and 
\cref{lemma:ineqforB}}\notag\\
&=\left(1-\frac u2\frac{(1+2u)^3}{(1+u)^2}\phi(u)-
\sqrt{1-u\frac{(1+2u)^2(1+3u)}{(1+u)}\phi(u)}\right)/u^2=:\theta(u),\\
|C|&\le\frac{\frac{(1-2u)(1+2u)^3}{2(1+u)}\phi(u)+(1+2u)(1+u)(2\sqrt{1+2u}-1)\theta(u)}{1-u\frac{(1-u)(1+2u)^4}{2(1+u)^2}{\phi(u)}}:=\chi(u),\label{ineq:C}\\
&\qquad\text{from 
\cref{eq:def-C,ineq:rho,ineq:thetatau,ineq:s-low,ineq:s-up,ineq:B}}\notag\\
|\kappa|
&\le \frac{\frac{1}{\frac{1+u}{1+2u}-u^2\chi(u)}-1-u}{u^2}-7=:\psi(u)
\qquad\text{from \cref{ineq:C} and \cref{thm-alg-Borges}}\label{ineq:R}.
\end{align}
\end{small}
The derivation of \cref{ineq:B} relies on the following lemma that captures
some
of the correlations.
\begin{lemma}Any real numbers $x,y$ such that $|x|+|y|\le1$ satisfy
the
inequality\label{lemma:ineqforB}
\[\left|1+x/2-\sqrt{1+x+y}\right|\le 1-|x|/2-\sqrt{1-|x|-|y|}.\]
\end{lemma}
\begin{proof} Let 
\[f(x,y)=1+x/2-\sqrt{1+x+y}.\]
We first show the inequality with $f$ in place of $|f|$.
Since the square root is increasing, we have
\[f(x,y)\le f(x,|y|),\qquad x-|y|\ge -1.\]
Next, if \(-1\le x\le 0\), then \(|x|=-|x|\) and the inequality follows.
Otherwise, \(|x|=x\) and the inequality follows from showing that
\(g(x)\ge0\) with
\begin{align*}
g(x)&=1-x/2-\sqrt{1-x-|y|}-\left(1+x/2-\sqrt{1+x-|y|}\right)\\
&=\sqrt{1+x-|y|}-\sqrt{1-x-|y|}-x.
\end{align*}
We have \(g(0)=0\) and the derivative is
\[g'(x)=\frac1{2\sqrt{1-x-|y|}}+\frac1{2\sqrt{1+x-|y|}}-1.\]
The value at 0 is
\[g'(0)=\frac1{\sqrt{1-|y|}}-1\ge0\]
and its derivative is
\[g''(x)=\frac{(1-x-|y|)^{-3/2}}{4}-\frac{(1+x-|y|)^{-3/2}}4,\]
which is \(\ge 0\) for \(0\le x\le1\) since the function
\(x\mapsto x^{-3/2}\) is decreasing. Therefore \(g'\ge0\) and \(g\ge0\).

For the absolute value, only the case when $f\le0$ remains to
be proved. In that case, the inequality to be proved becomes
\[-1-x/2+\sqrt{1+x+y}\le 1-|x|/2-\sqrt{1-|x|-|y|}.\]
From 
\[\sqrt{1+u}\le1+u/2,\qquad u\ge-1,\]
the left-hand side is upper bounded by $y/2$, while the right-hand
side is lower bounded by $|y|/2$, which concludes the proof.
\end{proof}
The application of this lemma in \cref{ineq:B} requires 
\[\left|u\frac{(1+2u)^3}{(1+u)^2}\phi
(u)\right|+\left|u^2\frac{(1+2u)^2(2+3u)}{(1+u)^2}\phi
(u)\right|\le1.\]
Since $\phi$ is positive for $0\le u\le1/4$, this is equivalent to
\[\left(u\frac{(1+2u)^3}{(1+u)^2}+u^2\frac{(1+2u)^2(2+3u)}{(1+u)^2}\right)\phi(u)\le 1.\]
The function on the left-hand side is strictly increasing and reaches
the value~1 at
$u\simeq0.1110831553$. Thus the inequality holds for $u\le
1/2^4=0.0625$.

\Cref{ineq:R} in the last line of the derivation above shows that 
the 
relative
error of Algorithm 4 is bounded by the function 
$u+(7+\kappa) u^2$
with
\[\kappa(u)=\frac{\frac{1}{\frac{1+u}{1+2 u}-\frac{u^{2} \left(-
\frac{\left
(1+2 u \right)^{3} \left(-1+2 u \right) \phi(u)}{2 \left(1+u
\right)}+\left(1+2 u \right) \left(1+u \right) \left(2 \sqrt{1+2
u}-1\right) \theta(u) \right)}{1+\frac{u \left(u -1\right) \left(1+2 u
\right)^{4} \phi(u)}{2 \left(1+u \right)^{2}}}}-1-u}{u^2}-7=O(u),\]
with $\phi$ and $\theta$ from \cref{ineq:thetatau,ineq:B}.
By analysing the univariate function $\kappa(u)$, for instance via
a polynomial that it cancels, we
find that it
is increasing for \(u\le 0.1\). Bounds are thus obtained by evaluating
numerically \(\kappa\) at the upper end of the interval of interest,
giving the values in the theorem.

\section{Proof of Theorem \ref{thm:error-Kahan}}
\label{sec:proof-kahan}
We detail the steps of the proof of \cref{thm:error-Kahan}, while omitting
expressions that
are too large to be reproduced here\footnote{A Maple session
detailing the computation is available on
arXiv with this article.}.

\subsection{Step-by-step analysis}

For this analysis, we assume that $u\le1/32$. One could accommodate larger values of $u$ by splitting the interval further, but this would only make the analysis more technical. 

When $x/2\le y\le x$, by \cref{lemma-sterbenz}, the
first step is
exact:
$$
\delta=x-y=x(1-\alpha),
$$
with $\alpha=y/x$ as in the other analyses and now $\alpha\in[1/2,1]$. Next 
$$
r_2=\frac\delta y(1+\epsilon_{r_2}u)\in[0,1],
$$
with $|\epsilon_{r_2}|\le 1-2u$ by \cref{lemma:jeannerod-rump}. 

Later on during the analysis, it is useful to split the interval $\alpha\in[1/2,1]$ into two subintervals $[1/2,3/5],[3/5,1]$. We use the exponents $(1),(2)$ to denote the values of the variables for $\alpha$ in each of these subintervals. For instance,
$$
r_2^{(2)}\in\left[0,\operatorname{RN}\left(\frac23\right)\right],\quad
r_2^{(1)}\in\left[\operatorname{RN}\left(\frac23\right),1\right],
$$
with $|\operatorname{RN}(2/3)-2/3|\le u/2$.

The next step gives $2r_2$ without error. Next, 
$$
r_3=(r_2^2+2r_2)(1+\epsilon_{r_3}u)\in[0,3],
$$
with $|\epsilon_{r_3}|\le 1/(1+u)$ by \cref{lemma-knuth-dekker}. Propagating intervals and
using absolute errors at the endpoints gives
$$
r_3^{(2)}\in[0,\left(\frac23+\frac u2\right)^2+\frac43+2u],\quad
r_3^{(1)}\in[\left(\frac23-\frac u2\right)^2+\frac43-2u,3]
.
$$
As $2+r_3\in[2,5]$, for the next step, we use an absolute error:
$$
r_4=2+r_3+{\epsilon_{r_4}}u,
$$
with $|\epsilon_{r_4}|\le 2$ for $r_4^{(2)}$ and $|\epsilon_{r_4}|\le 4$ for $r_4^{(1)}$. Thus,
$$
r_4^{(2)}\in[2,\left(\frac23+\frac u2\right)^2+\frac{4}3+2u],\quad
r_4^{(1)}\in[\left(\frac23-\frac u2\right)^2+\frac{10}3-4u,5].
$$
For the next step again, as $\sqrt{r_4}\in[\sqrt{2},\sqrt{5}]\subset[1,4]$, we use an absolute error
$$
s_2=\sqrt{r_4}+\epsilon_{s_2}u,
$$
with $|\epsilon_{s_2}|\le 1$ for $s_2^{(2)}$ while $|\epsilon_{s_2}|\le 2$  for $s_2^{(1)}$. At this stage,
$$
s_2^{(2)}\in[\sqrt2-u,2],\quad
s_2^{(1)}\in[\frac{\sqrt{136-168u+9u^2}}{6}-u,\sqrt5+2u].
$$
Next, we have $R_2=\sqrt2+\epsilon_{R_2}u$, with $|\epsilon_{R_2}|\le 1$ and similarly $P_h=1+\sqrt2+\epsilon_{P_h}u$ with $|\epsilon_{P_h}|\le 2$ and $P_\ell=1+\sqrt2-P_h+\epsilon_{P_\ell}u^2$ satisfies $|P_\ell|<2u$ and $|\epsilon_{P_\ell}|\le 1$.

Thus $R_2+s_2\in[2\sqrt2-2u,\sqrt2+\sqrt5+3u]\subset(2,4)$, so that 
$$
d=R_2+s_2+\epsilon_du,\qquad\text{with $|\epsilon_d|\le 2$}.
$$
Now,
$$
d^{(2)}\in[2\sqrt2-4u,\sqrt2+2+3u],\quad
d^{(1)}\in[\sqrt2+\frac{\sqrt{136-168u+9u^2}}{6}-4u,\sqrt2+\sqrt5+5u].
$$
For the next three steps, using relative errors gives
$$
q=\frac{r_3}d(1+\epsilon_qu),\quad
r_5=(q+P_\ell)(1+\epsilon_{r_5}u),\quad
r_6=(r_5+r_2)(1+\epsilon_{r_6}u),
$$
with $|\epsilon_q|\le 1-2u$ and $|\epsilon_{r_5}|,|\epsilon_{r_6}|$ both bounded by $1/(1+u)$. Propagating intervals shows that
\begin{align*}
&q^{(2)}\in[0,\frac{\left(\frac23+\frac u2\right)^2+\frac43+2u}{2\sqrt2-4u}+\frac u2],\\
&q^{(1)}\in[\frac{\left(\frac23-\frac u2\right)^2+\frac43-2u}{\sqrt2+\sqrt5+5u}-\frac u4,
\frac3{\sqrt2+\frac{\sqrt{136-168u+9u^2}}{6}-4u}+\frac u2],\\
&r_5^{(2)}\in[-2u,\frac{\left(\frac23+\frac u2\right)^2+\frac43+2u}{2\sqrt2-4u}+3u],\\
&r_5^{(1)}\in[\frac{\left(\frac23-\frac u2\right)^2+\frac43-2u}{\sqrt2+\sqrt5+5u}-\frac 52u,
\frac3{\sqrt2+\frac{\sqrt{136-168u+9u^2}}{6}-4u}+\frac 72u],\\
&r_6^{(2)}\in[-2u,\frac{\left(\frac23+\frac u2\right)^2+\frac43+2u}{2\sqrt2-4u}+\frac23+\frac92 u],\\
&r_6^{(1)}\in[1,1+\frac3{\sqrt2+\frac{\sqrt{136-168u+9u^2}}{6}-4u}+\frac {11}2u].
\end{align*}
It follows that 
$$
r_6^{(2)}+P_h\subset[2,4],\quad
r_6^{(1)}+P_h\subset[3,5],
$$
 so that one can use an absolute error
$$
z=r_6+P_h+\epsilon_zu,
$$
with $|\epsilon_z|\le 2$ for $z^{(2)}$ and $|\epsilon_z|\le 4$ for $z^{(1)}$.

For the last two steps, using relative errors gives
$$
z_2=\frac{y}z(1+\epsilon_{z_2}u),\quad \rho=(x+z_2)(1+\epsilon_\rho u),
$$
with $|\epsilon_{z_2}|\le 1-2u$ and $|\epsilon_\rho|\le1/(1+u).$ This concludes the step-by-step analysis.

\subsection{Partial derivatives of the relative error}

We now consider the relative error
$$
R=\frac\rho{\sqrt{x^2+y^2}}-1,
$$
where $\rho$ is the value computed by the algorithm and $R$ is viewed as a function of the variables $(x,\alpha,\epsilon_{r_2},\epsilon_{r_3},\epsilon_{r_4},\epsilon_{s_2},\epsilon_{R_2},\epsilon_d,\epsilon_q,\epsilon_{P_\ell},\epsilon_{P_h},\epsilon_{r_5},\epsilon_{r_6},\epsilon_z,\epsilon_{z_2},\epsilon_\rho).$

Remarkably, the design of the algorithm that avoids all subtractions leads to a simple untangling of many of the correlations between errors, summarized by the following lemma.

\begin{lemma} For $u\le1/32$, one has
\begin{gather*}
\frac{\partial R}{\partial \epsilon_{r_2}}\le 0,\quad
\frac{\partial R}{\partial \epsilon_{r_3}}\le 0,\quad
\frac{\partial R}{\partial \epsilon_{r_4}}\ge 0,\quad
\frac{\partial R}{\partial \epsilon_{s_2}}\ge 0,\quad
\frac{\partial R}{\partial \epsilon_{R_2}}\ge 0,\quad
\frac{\partial R}{\partial \epsilon_{d}}\ge 0,\quad
\frac{\partial R}{\partial \epsilon_{q}}\le 0,\\
\frac{\partial R}{\partial \epsilon_{P_\ell}}\le 0,\quad
{r_5}\frac{\partial R}{\partial \epsilon_{r_5}}\le 0,\quad
r_6\frac{\partial R}{\partial \epsilon_{r_6}}\le 0,\quad
\frac{\partial R}{\partial \epsilon_{z}}\le 0,\quad
\frac{\partial R}{\partial \epsilon_{z_2}}\ge 0,\quad
\frac{\partial R}{\partial \epsilon_{\rho}}\ge 0.
\end{gather*}
\end{lemma}
\begin{proof} For sign questions, it is sufficient to
consider the partial derivatives of $\rho$ rather than $R$, except for $\alpha$. The signs in the lemma follow from closed-forms expressions of the derivatives obtained by the chain rule:
\begin{align*}
\frac{\partial \rho}{\partial\epsilon_\rho}&=\frac{\rho u}{1+ u\epsilon_\rho},\\
\frac{\partial \rho}{\partial\epsilon_{z_2}}&=\frac{y u}z(1+u\epsilon_\rho),\\
\frac{\partial \rho}{\partial\epsilon_z}&=-\frac{y u}{z^2}(1+u\epsilon_{z_2})(1+u\epsilon_\rho),\\
\frac{\partial \rho}{\partial\epsilon_{r_6}}&=-\frac{yur_6}{z^2}\frac{(1+u\epsilon_{z_2})(1+u\epsilon_\rho)}{1+u\epsilon_{r_6}},\\
\frac{\partial\rho}{\partial\epsilon_{r_5}}&=-\frac{yur_5}{z^2}\frac{(1+u\epsilon_{r_6})(1+u\epsilon_{z_2})(1+u\epsilon_\rho)}{1+u\epsilon_{r_5}},\\
\frac{\partial\rho}{\partial\epsilon_{P_\ell}}&=-\frac{yu^2}{z^2}(1+u\epsilon_{r_6})({1+u\epsilon_{r_5}})(1+u\epsilon_{z_2})(1+u\epsilon_\rho),\\
\frac{\partial\rho}{\partial\epsilon_{q}}&=-\frac{yur_3}{z^2d}(1+u\epsilon_{r_5})(1+u\epsilon_{r_6})(1+u\epsilon_{z_2})(1+u\epsilon_\rho),\\
\frac{\partial\rho}{\partial\epsilon_{d}}&=
\frac{\partial\rho}{\partial\epsilon_{s_2}}=\frac{\partial\rho}{\partial\epsilon_{R_2}}=
\frac{yur_3}{z^2d^2}(1+u\epsilon_q)(1+u\epsilon_{r_5})(1+u\epsilon_{r_6})(1+u\epsilon_{z_2})(1+u\epsilon_\rho),\\
\frac{\partial\rho}{\partial\epsilon_{r_4}}&=\frac{yur_3}{2z^2d^2\sqrt{r_4}}(1+u\epsilon_q)(1+u\epsilon_{r_5})(1+u\epsilon_{r_6})(1+u\epsilon_{z_2})(1+u\epsilon_\rho),\\
\frac{\partial\rho}{\partial\epsilon_{r_3}}&=-\frac{yur_3}{z^2d(1+u\epsilon_{r_3})}\left(1-\frac{r_3}{2d\sqrt{r_4}}\right)(1+u\epsilon_q)(1+u\epsilon_{r_5})(1+u\epsilon_{r_6})(1+u\epsilon_{z_2})(1+u\epsilon_\rho),\\
\frac{\partial\rho}{\partial\epsilon_{r_2}}&=-\frac{\delta u}{z^2}
\left(1+\frac{2(1+r_2)}{d}\left(1-\frac{r_3}{2d\sqrt{r_4}}\right)(1+u\epsilon_q)(1+u\epsilon_{r_3})(1+u\epsilon_{r_5})\right)
(1+u\epsilon_{r_6})(1+u\epsilon_{z_2})(1+u\epsilon_\rho).
\end{align*}
The only sign that is not directly obvious from the previous identities is that of the factor
\begin{align*}
1-\frac{r_3}{2d\sqrt{r_4}}=1-\frac{r_3}{2d\sqrt{2+r_3+u\epsilon_{r_4}}}&\ge
1-\frac{r_3}{2(2\sqrt{2}-4u)\sqrt{2+r_3-4u}}\\
&\ge 1-\frac{3}{2(2\sqrt{2}-4u)\sqrt{5-4u}}\ge0.\qedhere
\end{align*}
\end{proof}

\subsection{Discussion on \texorpdfstring{$r_5,r_6$}{r5,r6}}
As \(r_2\ge0\) and the sign of \(r_6\) is that of \(r_2+r_5\), there are
only three possibilities for the signs of \(r_5\) and \(r_6\):
\((+,+),(-,-),(-,+)\). The maximal value of \(R\) is upper bounded by
the value it reaches at one of the following corresponding values of
\((\epsilon_{r_5},\epsilon_{r_6})\):
\[\phi_{++}=(-\frac1{1+u},-\frac1{1+u}),\quad
\phi_{--}=(\frac1{1+u},\frac1{1+u}),\quad
\phi_{-+}=(\frac1{1+u},-\frac1{1+u}).\]
Similarly, the minimal value of \(R\) is lower bounded by the value it
reaches at the opposite of these points.

A computation similar to that of the derivatives above gives the
derivative:
\begin{align*}
\frac{\partial\rho}{\partial\epsilon_{P_h}}&=\frac{yu^2}{z^2}\left(\epsilon_{r_5}+\epsilon_{r_6}+u\epsilon_{r_5}\epsilon_{r_6}\right)(1+u\epsilon_{z_2})(1+u\epsilon_\rho).
\end{align*}
so that the sign of this derivative at the extremal points above are
easily obtained and thus also the corresponding extremal values of
\(\epsilon_{P_h}\). This leads to three possible points to be considered
for \((\epsilon_{r_5},\epsilon_{r_6},\epsilon_{P_h})\):
\begin{align}
\psi_{++}&=(-\frac1{1+u},-\frac1{1+u},-2),\label{eq:phipp}\\ 
\psi_{--}&=(\frac1{1+u},\frac1{1+u},2),\label{eq:phimm}\\
\psi_{-+}&=(\frac1{1+u},-\frac1{1+u},-2).\label{eq:phimp}
\end{align}
Their opposite is to be used for the minimal value of \(R\).

\subsection{Extremal values}
As a consequence of the previous discussion, the maximal value of \(R\)
is upper bounded by the value it reaches at
\begin{multline*}
(\epsilon_{r_2},\epsilon_{r_3},\epsilon_
{r_4},\epsilon_{s_2},\epsilon_{R_2},\epsilon_d,\epsilon_q,\epsilon_{P_\ell},\epsilon_z,\epsilon_{z_2},\epsilon_\rho)=\pi_+,\quad (\epsilon_{r_5},\epsilon_{r_6},\epsilon_{P_h})\in\{\psi_{++},\psi_{--},\psi_{-+}\}\\
\pi_+:=\left(-1+2u,-\frac1{1+u},e,\frac e2,1,2,-1+2u,-1,-e,1-2u,\frac1{1+u}\right),
\end{multline*}
with \((\psi_{++},\psi_{--},\psi_{-+})\) from 
\cref{eq:phipp,eq:phimm,eq:phimp}, \(e=2\) if
\(\alpha\ge3/5\) and \(e=4\) otherwise. Similarly, its minimal value is
lower bounded by the value it reaches at \(\pi_-=-\pi_+\) and one of
\(\{-\psi_{++},-\psi_{--},-\psi_{-+}\}\). The problem is thus reduced to
functions of only the two variables \((\alpha,u)\). In each case, the
sign of the derivative with respect to \(\alpha\) is easily predicted by
plotting its graph with high precision for
\((\alpha,u)\in[3/5,1]\times[0,1/32]\). A computationally intensive part
of this analysis is to actually prove the following.

\begin{lemma} At each of the three points
\((\pi_+,\psi_{++}),(\pi_+,\psi_{--}),(\pi_+,\psi_{-+})\), for
\((\alpha,u)\in[3/5,1]\times[0,1/32]\), the derivative of \(R\) wrt
\(\alpha\) is positive.
\end{lemma}
\begin{proof} First, a computation similar to that of the other
derivatives gives the following expression for
\(S:=\frac1{R+1}\frac{\partial R}{\partial\alpha}\), that has the same
sign as \(\partial R/\partial\alpha\):
\begin{multline*}
\frac{
\Bigl(1+\frac{\Bigl(1+\frac{2 \left(r_{2}+1\right) 
\left(1-\frac{r_{3}}{2d \sqrt{r_{3}+2+u \epsilon_{r_{4}}}}\right)
\left(1+u \epsilon_{r_{3}}\right)
\left(1+u \epsilon_{q}\right) 
\left(1+u \epsilon_{r_{5}}\right)}
 {d}\Bigr) \left(1+u \epsilon_{r_{2}}\right) 
 \left(1+u \epsilon_{r_{6}}\right)}{\alpha  z}\Bigr) \left(1+u
 \epsilon_{z_{2}}\right)}{z \left(1 +z_{2}/x\right)}\\
-\frac{\alpha}{\alpha^{2}+1}.
\end{multline*}
Since \(\alpha\ge3/5\), we use \(\pi_+\) with \(e=2\). At this point,
and for each of the three choices for \(\psi\), this expression gives us
a function \(S(\alpha,u)\). That function is a root of a (large)
polynomial of degree only 2
\[A(\alpha,u)S^2+B(\alpha,u)S+C(\alpha,u)=0\]
the polynomials \(A,B,C\) having their coefficients in
\(\mathbb Q(\sqrt2)\), and degree 16 in \(\alpha\) and between \(32\)
and \(36\) in \(u,\) depending on the three cases. All three cases can
be dealt with in the same way. First, by direct computation, one sees
that the value of \(S\) at \((\alpha,u)=(1,1/32)\) is positive. Next, if
\(S\) becomes 0 in the rectangle \((\alpha,u)\in[3/5,1]\times[0,1/32]\),
one has that \(C\) becomes 0 there as well. So it is sufficient to prove
that this does not happen. This replaces proving the positivity of an
algebraic function to that of a polynomial in the same number of
variables. On each of the four sides of the rectangle, the polynomial
\(C\) specializes to a univariate polynomial whose positivity is easy to
check. Next, one considers the polynomial system
\(\{C,\partial C/\partial u,\partial C/\partial\alpha\}\) which is
satisfied at the extrema of \(C\). For the cases \(\psi_{--}\) and
\(\psi_{-+}\), it turns out to be sufficient to compute the resultant of
the last two polynomials with respect to \(\alpha\). This can be
achieved in less than 10 sec. and yields a univariate polynomial (of
degree 569 in the first case and 524 in the second one) that can then be
seen not to have a zero for \(u\in(0,1/32)\). In the case of
\(\psi_{++}\), there are zeroes. More information is obtained by
computing a parameterization of the roots of the system
\(\{s^2-2,C,\partial C/\partial u,\partial C/\partial\alpha\}\), where
\(\sqrt2\) is replaced by \(s\) in the polynomials so that they all lie
in \(\mathbb Q[s,\alpha,u]\). This parameterization is computed by the
\texttt{msolve} library~\cite{BerthomieuEderDin2021} in 
40~seconds. It has the
form
\[P(u)=0,\quad P'(u)\alpha+Q_\alpha(u)=0,\quad P'(u)s+Q_s(u)=0,\]
with \(P\) a polynomial of degree 530, while \(Q_\alpha\) and \(Q_s\)
have degree 529. The root of the resultant above is still a root of
\(P\). It belongs to the interval \([0.00443,0.00445]\). Evaluating
\(-Q_\alpha/P'\) at this point by interval arithmetic shows that the
corresponding value of \(\alpha\) is negative. This
shows that there are no extrema of \(C\) inside the rectangle and
therefore its minimal value is reached on the boundary where it is
positive, showing that \(S\) does not vanish in the rectangle and
establishing its positivity.
\end{proof}
\begin{corollary}For \((\alpha,u)\in[3/5,1]\times[0,1/32]\), the
relative error \(R\) is upper bounded by $\phi(u)$ from \cref{thm:error-Kahan}.
\end{corollary}
\begin{proof} The consequence of the previous lemma is that the
relative error \(R\) is bounded by its values at \(\pi_+\), each of
\((\psi_{++},\psi_{--},\psi_{-+})\) and \(\alpha=1\). This gives three
rational functions in only the variable \(u\) to compare. A direct
computation then shows that the maximal value is reached at
\(\psi_{--}\), where the bound is given by $\phi(u)$ from \cref{thm:error-Kahan}.
\end{proof}

\subsection{Final steps in the proof of Theorem \ref{thm:error-Kahan}}
To conclude the proof for \(\alpha\ge3/5\), we also need to check that
the value of \(R\) at \(-\pi_+\), does not exceed \(\phi(u)\) in
absolute value; that the bound for $\alpha\le 3/5$ is smaller
than~$\phi(u)$ and similarly for the other path of the algorithm when
$\delta>y$. 

There is no new difficulty there, and the computation
follows the same steps, so as to control the size of intermediate
expressions when more than two variables are involved. 

\end{document}